# A Classification of Flexible Kokotsakis Polyhedra with Reducible Quadrilaterals

YANG LIU, Academy of Mathematics and System Science, CAS, China

We study a class of mechanisms known as Kokotsakis polyhedra with a quadrangular base. These are $3 \times 3$ quadrilateral meshes whose faces are rigid bodies and joined by hinges at the common edges. In contrast to existing work, the quadrilateral faces do not necessarily have to be planar. In general, such a mesh is rigid. The problem of finding and classifying the flexible ones is old, but until now largely unsolved. It appears that the tangent values of the dihedral angles between different faces are algebraically related through polynomials. Specifically, this article deals with the case when these polynomials are reducible. We explore the conditions for reducibility to characterize all possible shape restrictions that lead to flexible polyhedra.

## 1 INTRODUCTION

In many disciplines within the sciences and engineering, there is a growing interest in flexible or deployable structures. These provide a continuous morphing between different states of a shape and found widespread utility from robotics and solar cells to meta-materials and architecture.

Deployable or flexible structures may contain non-rigid parts. However, in this paper we are only dealing with mechanisms in the classical sense. They are formed by rigid bodies connected by hinges. Our interest is in flexible quadrilateral meshes, whose faces are rigid and joined by hinges at the common edges. It is well-known that such a quadrilateral mesh is a mechanism if and only if any $3 \times 3$ sub-mesh is flexible (a $2 \times 2$ mesh is always flexible). This has been shown by Schief et al. in Theorem 3.2 of [9], who formulate the problem for planar faces, but planarity does not enter their proof. Therefore, the first key step in the determination of all flexible quadrilateral meshes is the *classification of all flexible* $3 \times 3$ *quadrilateral meshes*, which is the main contribution of the present paper.

The problem of flexible polyhedra can be traced back to the last century. Very similar mechanisms to the ones studied in our paper have been introduced by Kokotsakis [4]. His mechanisms possess a central planar, but not necessarily quadrilateral face, surrounded by a belt of planar faces. The work of Kokotsakis inspired further research on such structures [12], in particular with a focus on origami [13–15].

A breakthrough in this area has been made by Izmestiev [2] who derived a complete classification of flexible $3 \times 3$ quadrilateral meshes with planar faces. It is, however, still unknown which of these types can be parts of larger flexible $n \times m$ quadrilateral meshes. Larger flexible quadrilateral meshes with planar faces have been studied by Schief et al. [9] as discrete integrable systems and counterparts to isometric deformations of conjugate nets on smooth surfaces. Their results are largely rooted in second-order infinitesimal flexibility. A wealth of results on first-order infinitesimal flexibility is found in the book by Sauer [7]. There, we also find a detailed study of special quadrilateral meshes with planar faces that are mechanisms. They come in two types. One type, the so-called Voss meshes, are discrete counterparts to surfaces with a conjugate net of geodesics, first described by Voss [16] (see also [9]). The other type is the T-nets, first described by Graf and Sauer [8] (for a recent description in English see [11]), are discrete counterparts to an affine generalization of moulding surfaces. A very special case of T-nets is the famous Miura origami structures [5]. T-nets are also capable of forming tubular flexible structures and flat-foldable metamaterials [10]. The isometric deformations of T-nets and their smooth counterparts have recently been studied by Izmestiev et al. [3].

Author's address: Yang Liu, Academy of Mathematics and System Science, CAS, Beijing, China, 100000, yangrio@amss.ac.cn.

The question of the existence of a flexible Kokotsakis polyhedron with skew-quadrilateral faces has been posed by Sauer [7], but remained unsolved for a long time. The first non-trivial example of such a polyhedron was constructed by Nawratil in 2022 [6], where he built several flexible meshes of 'isogonal type'. That means some specific shape constraints were added to their skew faces. However, the contribution of [6] only provides a few particular solutions without a formal analysis of isogonal meshes.

In this paper, we are going to complete Nawratil's work [6] by giving a thorough analysis of the isogonal type and providing a systematic construction that includes all possible flexible isogonal meshes. Beyond this, we will also introduce new flexible members, the singular type, and the corresponding construction. This work will be a great starter for the complete classification on general flexible $3 \times 3$ meshes with skew faces.

### 1.1 Main approach and contributions

The main mechanism we study in this paper is the flexible mesh[1] in Fig.1 left, which contains nine quadrilaterals linked by hinges that allow for the linked faces to rotate. A given mesh uniquely determines a system whose solution set is the trajectory of the motion of the mesh. It will be shown that the trajectory is generally a curve in 4-dimensional space described by a polynomial system

$$M : \{G^{(1)}(x_1, x_2) = G^{(2)}(x_2, x_3) = G^{(3)}(x_3, x_4) = G^{(4)}(x_4, x_1) = 0\} \tag{1}$$

where the explicit form of $G^{(i)}$ will be given in Definition 2.4. We split the system into two subsystems $S_1 : \{G^{(1)} = G^{(2)} = 0\}$ and $S_2 : \{G^{(3)} = G^{(4)} = 0\}$, hence the solution set of $M$ equals the following fiber product

$$\text{Sol.}(M) = \text{Sol.}(S_1) \underset{\{x_1,x_3\}}{\times} \text{Sol.}(S_2) := \{(x_i)_{i=1}^4 | (x_1, x_2, x_3) \in \text{Sol.}(S_1), (x_3, x_4, x_1) \in \text{Sol.}(S_2)\}. \tag{2}$$

Generally, four equations of four unknowns lead to finitely many solutions. For a flexible mesh, however, Sol.$(M)$ contains infinitely many solutions, which requires the $(x_1, x_3)$-projections of Sol.$(S_1)$ and Sol.$(S_2)$ overlap (not just intersect). In particular, if we are searching Sol.$(M)$ in $\mathbb{C}^4$, it is equivalent to saying that the following resultants share a common factor, i.e.

$$\gcd(\text{Res}(G^{(1)}, G^{(2)}; x_2), \text{Res}(G^{(3)}, G^{(4)}; x_4)) \neq 1.$$

In this regard, we firstly study $S_i$ based on the irreducible components (algebraic varieties) of Sol.$(S_i)$. Afterwards, we discuss all possible combinations of $S_1$ and $S_2$ such that the fiber product in equation (2) is an infinite set. Finally, we obtain a classification and construction of flexible meshes through $M$.

Our main contribution is to the meshes in which both $S_i$ consist of reducible polynomials. That is, all $G^{(i)}$ are reducible. We study $S_i$ with such a property, and then provide constructions for all possible combinations of $S_1$ and $S_2$ to form a system $M$ which corresponds to a flexible mesh. However, in Section 6, we construct a so-called **constant mesh** that may contain irreducible $G^{(i)}$. This particular type is much easier to construct than the others.

PLEASE NOTE: Many proofs and arguments of this article use symbolic computation. We provide the readers a Maple™ script[2] to verify all nontrivial calculation results.

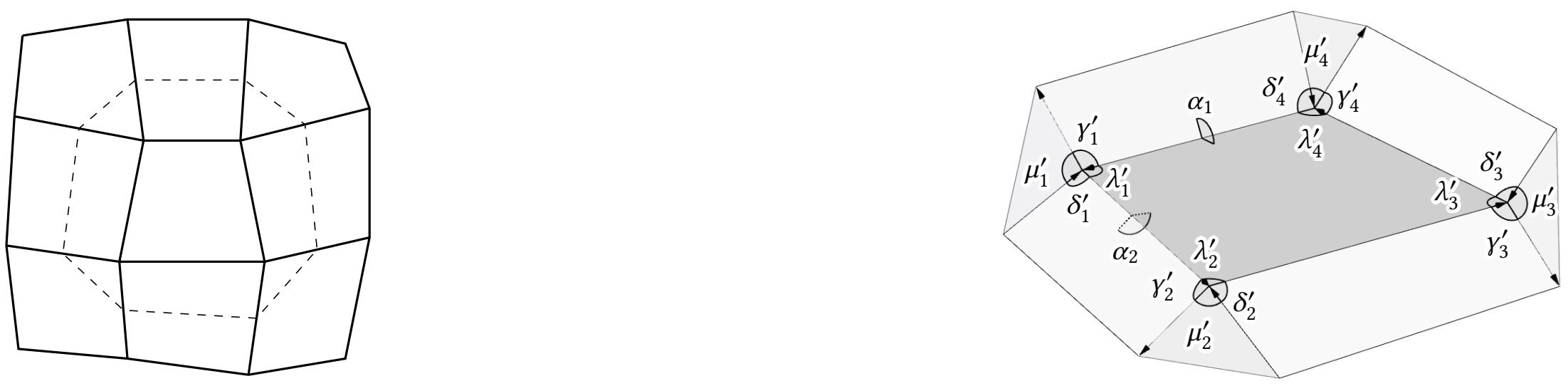

Fig. 1. Left: Sketch of a 3 by 3 quadrilateral mesh, or equivalently, Kokotsakis polyhedron with a quadrangular base; Right: A flexible mesh with fixed interior angles $\lambda'_i, \gamma'_i, \mu'_i, \delta'_i$ of planar faces. Flexible angles $\alpha_i$ are dihedral angles between the central face and the attached faces. The rigorous definitions of these angles are in Appendix A.

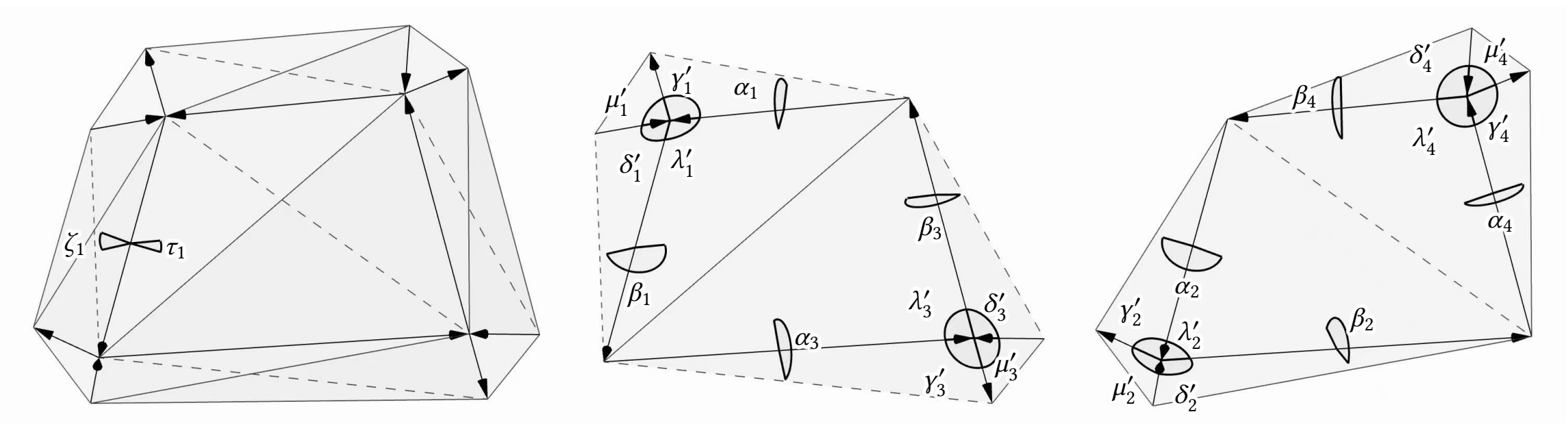

Fig. 2. A skew flexible mesh (left) and its decomposition (middle and right). $(\lambda'_i, \gamma'_i, \mu'_i, \delta'_i)$ are fixed angles among oriented edges; $\tau_i, \zeta_i$ are fixed dihedral angles along each common edge of the central tetrahedron and the attached tetrahedron respectively; $\alpha_i, \beta_i$ are flexible dihedral angles (in planar case $\alpha_{i+1} = \beta_i$). All angles can be defined in a similar way (see Appendix A).

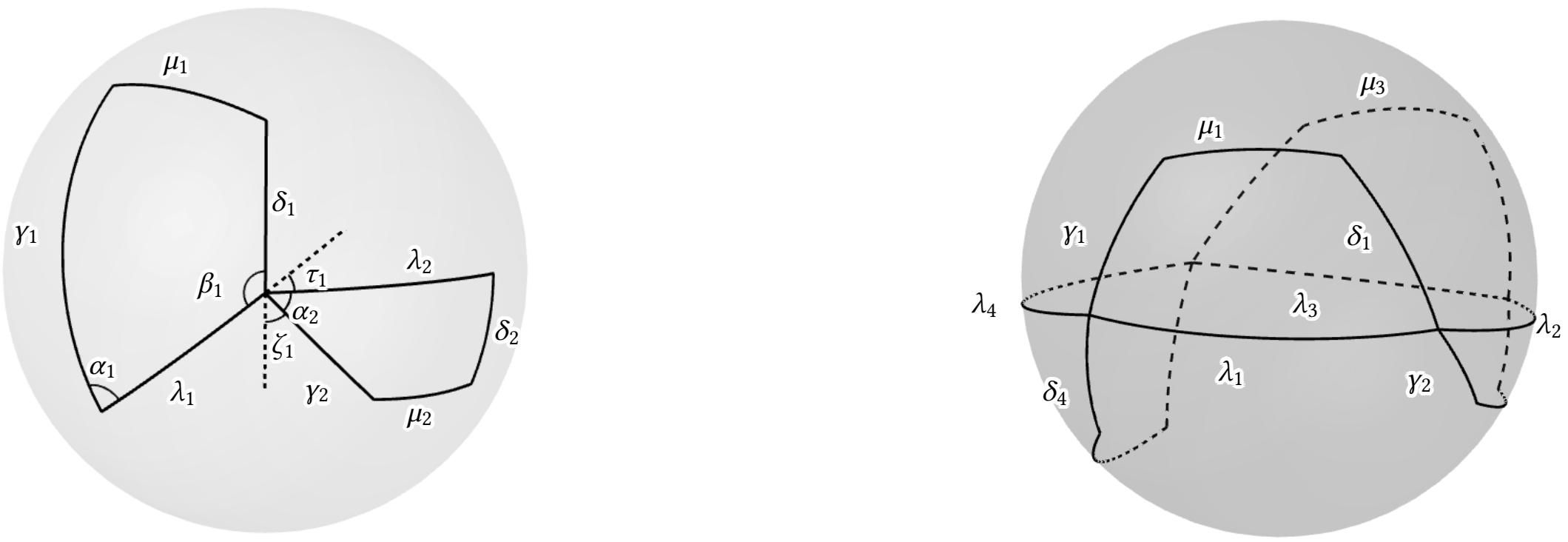

Fig. 3. Partial (left) and whole (right) spherical linkage of the mesh in Fig.2, $(\lambda_i, \gamma_i, \mu_i, \delta_i)$ and $(\lambda'_i, \gamma'_i, \mu'_i, \delta'_i)$ are complementary to $\pi$ respectively, the gap between $\beta_1$ and $\alpha_2$ is caused by $\tau_1 + \zeta_1$.

## 2 BASIC SETUP

### 2.1 From $\mathbb{R}^3$ to sphere

From Fig.1 left, it is obvious that the flexibility of a mesh only depends on the structure of a neighborhood of the central quadrilateral. So, without loss of generality, we can trim the mesh such that its corner faces are just triangles. Fig.1 right is an illustration of a flexible mesh with 9 planar quadrilaterals. Similarly, if the mesh contains nine skew quadrilaterals, each quadrilateral can be treated as a rigid tetrahedron (a skew quadrilateral becomes a tetrahedron once diagonals are added). Hence Fig.1 right can be generalized to Fig.2 left of a $3 \times 3$ tetrahedron mesh.

The classic approach to analyze the flexibility is to consider a spherical quadrilateral at each vertex of the central face, which is defined by oriented edges: In Fig.2, by collecting all vectors at the origin, the directions (of vectors) intersect the unit sphere and determine a **spherical linkage** of the mesh, see Fig.3 right, which contains 4 spherical quads.[3] Every 'edge' of each quad is an arc of a great circle. Particularly, each quad of the lingkage, say $(\lambda_1, \gamma_1, \mu_1, \delta_1)$, is associated with the corner surrounded by $(\lambda_1', \gamma_1', \mu_1', \delta_1')$ in Fig.2. It is easy to see that $(\lambda_i, \gamma_i, \mu_i, \delta_i)$ and $(\lambda_i', \gamma_i', \mu_i', \delta_i')$ are complementary to $\pi$ respectively. Angles $\tau_i, \zeta_i$ are (interior) dihedral angles in the corresponding tetrahedra. All these angles will be rigorously defined in Appendix A. In summary, the mesh vs. its spherical linkage, angles between edges are transformed to arcs of the quads, angles between planes are transformed to angles between arcs, flexible hinges of adjacent quadrilaterals are transformed to flexible joints of adjacent arcs. Hence the flexibility of the mesh is equivalent to the flexibility of the spherical linkage with the deformation rules:

- All arc lengths of each quad are fixed;
- All angles $\tau_i$ between $\lambda_i$ and $\lambda_{i+1}$ are fixed;
- All angles $\zeta_i$ between $\delta_i$ and $\gamma_{i+1}$ are fixed.

### 2.2 From sphere to polynomials

According to Fig.3 left, by setting

$$x_i = \tan\left(\frac{\alpha_i}{2}\right),\ y_i = \tan\left(\frac{\beta_i}{2}\right),\ f_i = \tan\left(\frac{\zeta_i + \tau_i}{4}\right) \tag{3}$$

we have an invertible relation between $y_i$ and $x_{i+1}$: $y_i = \frac{(1-f_i^2)x_{i+1}+2f_i}{1-f_i^2-2f_i x_{i+1}}$, or

$$h^{(i)}(y_i, x_{i+1}) = (1 - f_i^2)(y_i - x_{i+1}) - 2f_i(y_i x_{i+1} + 1) = 0,\ f_i \in (-1, 1],$$

where $i = 1, 2, 3, 4$ and $x_5 = x_1$.

When $\lambda_i, \gamma_i, \mu_i, \delta_i \in (0, \pi)$, the arc lengths and angles of the quad in Fig.3 left satisfy a polynomial equation during its deformation. That means when the quad changes its shape while all arcs are attaching the sphere and remaining the same lengths, the interior angles $\alpha_i, \beta_i$ must obey the so-called Bricard equation [1] (the proof can be also found in [12]).

$$\begin{cases} A_i x_i^2 y_i^2 + B_i x_i^2 + C_i y_i^2 + D_i x_i y_i + E_i = 0, \text{ or} \\ g^{(i)}(x_i, y_i) := a_i x_i^2 y_i^2 + b_i x_i^2 + c_i y_i^2 + x_i y_i + e_i = 0 \end{cases} \tag{4}$$

[1] In this article, 'mesh' is always meant by default a $3 \times 3$ mesh with quadrilateral faces which are not necessarily planar.
[2] https://drive.google.com/file/d/1oYnGvNC10R-36iqErdaAw6CPUxRd3BBu/view?usp=sharing
[3] To distinguish the Euclidean quadrilateral, we will call spherical quadrilateral just quad.

where

$$\begin{cases} A_i = \cos(\lambda_i + \gamma_i + \delta_i) - \cos(\mu_i), \\ B_i = \cos(\lambda_i + \gamma_i - \delta_i) - \cos(\mu_i), \\ C_i = \cos(\lambda_i - \gamma_i + \delta_i) - \cos(\mu_i), \\ D_i = 4\sin(\gamma_i)\sin(\delta_i), \\ E_i = \cos(\lambda_i - \gamma_i - \delta_i) - \cos(\mu_i), \\ (a_i, b_i, c_i, e_i) = \left(\frac{A_i}{D_i}, \frac{B_i}{D_i}, \frac{C_i}{D_i}, \frac{E_i}{D_i}\right). \end{cases}$$

Now all flexible angles $\alpha_i, \beta_i$ (as known as $x_i, y_i$) are related by eight polynomial equations $\{h^{(i)} = g^{(i)} = 0\}_{i=1}^4$. However, $y_i$ becomes redundant after the elimination (by resultant)

$$\{G^{(i)}(x_i, x_{i+1}) := \operatorname{Res}(g^{(i)}, h^{(i)}; y_i) = 0\}_{i=1}^4 \tag{5}$$

And the above $G^{(i)}$ is what we have mentioned in system (1).

### 2.3 Shape restrictions of the quads

Obviously, a single quad $Q_i$ is in general flexible. Its shape is able to continuously change while all arc lengths $\lambda_i, \gamma_i, \mu_i, \delta_i$ are fixed. Here we do not consider the case where 0 or $\pi$ appears to be one of the lengths, which means the corresponding quadrilateral in Fig.2 degenerates to a triangle. It is therefore reasonable to have

- General assumption: $\boxed{\lambda_i, \delta_i, \mu_i, \gamma_i \in (0, \pi), \forall i = 1, 2, 3, 4}$

since they can be defined by inner product. Bricard equation (4) is therefore applicable to all spherical linkages of meshes.

We are particularly interested in two types of spherical quads in this article. Here we follow the definitions from [2].

*Definition 2.1.* Quad $Q_i$ with arc lengths $(\lambda_i, \delta_i, \mu_i, \gamma_i)$ is called **(anti)isogonal** or an **(anti)isogram** if it satisfies one of the corresponding conditions:

(i) $\lambda_i = \mu_i, \gamma_i = \delta_i$ (isogram); (ii) $\lambda_i + \mu_i = \gamma_i + \delta_i = \pi$ (antiisogram).

*Definition 2.2.* Quad $Q_i$ with arc lengths $(\lambda_i, \delta_i, \mu_i, \gamma_i)$ is called **(anti)deltoidal** or an **(anti)deltoid** if it satisfies one of the corresponding conditions:

(iii) $\lambda_i = \gamma_i, \mu_i = \delta_i$ (deltoid); (iv) $\lambda_i + \gamma_i = \mu_i + \delta_i = \pi$ (antideltoid);
(v) $\lambda_i = \delta_i, \gamma_i = \mu_i$ (deltoid); (vi) $\lambda_i + \delta_i = \gamma_i + \mu_i = \pi$ (antideltoid).

Simple calculations can show the following.

Lemma 2.3. *In Bricard equation* (4)*, the conditions for (anti)isograms and (anti)deltoids are respectively equivalent to:*

(i) $\Leftrightarrow b_i = c_i = 0$; *(ii)* $\Leftrightarrow a_i = e_i = 0$; *(iii)* $\Leftrightarrow c_i = e_i = 0$;
(iv) $\Leftrightarrow a_i = b_i = 0$; *(v)* $\Leftrightarrow b_i = e_i = 0$; *(vi)* $\Leftrightarrow a_i = c_i = 0$.

### 2.4 Reformulation in algebra

We have seen how a given mesh (Fig.1 left) is transformed step by step to a polynomial system (5). It is natural to conduct the flexibility analysis of the mesh directly on its corresponding polynomials due to the fact: the mesh is flexible $\Leftrightarrow$ system (5) has infinitely many solutions.[4] Thus, we need to reformulate the flexibility problem and reclaim the goal of this paper in an algebraic fashion.

The following notations are the most important ones that must be specifically defined and stick to their meanings throughout the article.

*Definition 2.4.* $g^{(i)} \in \mathbb{R}[x_i, y_i]$ is a polynomial in the form

$$g^{(i)}(x_i, y_i) := a_i x_i^2 y_i^2 + b_i x_i^2 + c_i y_i^2 + x_i y_i + e_i$$

where the coefficients satisfy

$$(1 - 4a_i e_i - 4b_i c_i)^2 > 64 a_i b_i c_i e_i, \ \ |b_i + c_i - a_i - e_i| < 1. \tag{6}$$

$h^{(i)} \in \mathbb{R}[y_i, x_{i+1}]$ is a polynomial determined by a constant $f_i \in (-1, 1]$,

$$h^{(i)}(y_i, x_{i+1}) := (1 - f_i^2)(y_i - x_{i+1}) - 2f_i(y_i x_{i+1} + 1).$$

Finally, $G^{(i)} \in \mathbb{R}[x_i, x_{i+1}]$ is the resultant

$$G^{(i)}(x_i, x_{i+1}) := \operatorname{Res}(g^{(i)}, h^{(i)}; y_i).$$

The indices are always considered modulo 4, e.g., $x_5 = x_1$.

Corollary 2.5. *Every $g^{(i)}$ satisfies $(b_i - a_i)(c_i - e_i) < \frac{1}{4}$ and $(c_i - a_i)(b_i - e_i) < \frac{1}{4}$.*

Proof. Notice that $|x + y| < 1 \Rightarrow xy < \frac{1}{4}$, and Definition 2.4 requires

$$|b_i + c_i - a_i - e_i| = |(b_i - a_i) + (c_i - e_i)| = |(c_i - a_i) + (b_i - e_i)| < 1.$$

□

**Remark 1** (Recover the mesh from coefficients)**.** Inequalities (6) are induced by equation (4) in which $0 < \lambda_i, \gamma_i, \mu_i, \delta_i < \pi$. Symbolic computation can show

$$(1 - 4a_i e_i - 4b_i c_i)^2 > 64 a_i b_i c_i e_i \Leftrightarrow (\sin(\lambda_i)\sin(\mu_i)\csc(\delta_i)\csc(\gamma_i))^2 > 0,$$
$$|b_i + c_i - a_i - e_i| < 1 \Leftrightarrow |\cos(\lambda_i)| < 1.$$

One can easily verify that (please note $A_i = a_i D_i = 4a_i \sin(\gamma_i)\sin(\delta_i)$)

$$\begin{cases} \lambda_i = \cos^{-1}(b_i + c_i - a_i - e_i), \\ \gamma_i = \cot^{-1}\left(\frac{b_i - c_i - a_i + e_i}{\sin(\lambda_i)}\right), \\ \delta_i = \cot^{-1}\left(\frac{c_i + e_i - a_i - b_i}{\sin(\lambda_i)}\right), \\ \mu_i = \cos^{-1}(\cos(\lambda_i + \gamma_i + \delta_i) - A_i). \end{cases}$$

It is obvious that, except for $\mu_i$, the above angles always have real values. However, symbolic computation can show

$$\sin^2(\mu_i) = \frac{[(1 - 4a_i e_i - 4b_i c_i)^2 - 64 a_i b_i c_i e_i][1 - (b_i + c_i - a_i - e_i)^2]}{[1 - 4(b_i - a_i)(c_i - e_i)][1 - 4(c_i - a_i)(b_i - e_i)]} > 0$$

by Corollary 2.5, so $\mu_i$ is also real.

On the other hand, it is well-known that every skew quadrilateral $ABCD$ is determined by six free parameters. However, since flexibility is invariant under similarity, we may always assume $|AB| = 1$. So it is enough for us to construct the central face of a mesh with five parameters $(\lambda_1', \lambda_2', \lambda_3', \lambda_4', \tau_1)$ where $\lambda_i' = \pi - \lambda_i$ (Fig. 2), after which $(\tau_2, \tau_3, \tau_4)$ are known. The remaining angles can be obtained from $\zeta_i = 4\tan^{-1}(f_i) - \tau_i$.

[4]We accept the flexibility even for complex solutions.

In general, there is no guarantee that $g^{(i)} = 0$ provides a real curve in $\mathbb{R}^2$, e.g. when $\mu_i = \lambda_i + \gamma_i + \delta_i$. Thus, we will develop our theory disregarding reality and embeddability into Euclidean space.

For convenience, we use ideal and zero set instead of system and solution set.

*Definition 2.6.* An ideal in the form

$$S = (g^{(i)}, h^{(i)}, g^{(i+1)}, h^{(i+1)}) \subset \mathbb{C}[x_i, x_{i+1}, x_{i+2}, y_i, y_{i+1}]$$

is called a **coupling**. An ideal in the form

$$M = (G^{(1)}(x_1, x_2), G^{(2)}(x_2, x_3), G^{(3)}(x_3, x_4), G^{(4)}(x_4, x_1)) \subset \mathbb{C}[x_1, x_2, x_3, x_4]$$

is called a **matching** if the zero set

$$Z(M) := \{(x_1, x_2, x_3, x_4) \in (\mathbb{C}P^1)^4 : f(x_1, x_2, x_3, x_4) = 0, \forall f \in M\}.^5$$

is an infinite set (please note $(\mathbb{C}P^1)^n \neq \mathbb{C}P^n$).

For $S_1 = (g^{(1)}, h^{(1)}, g^{(2)}, h^{(2)})$, it is clear that $G^{(1)}, G^{(2)} \in S_1$ hence the solution set of $\{G^{(1)} = G^{(2)} = 0\}$ coincides with the $(x_1, x_2, x_3)$-projection of $Z(S_1)$, which calls back to the subsystems we mentioned in Section 1.1. In summary, it has been shown that every mesh uniquely determines an ideal $M = (G^{(1)}, G^{(2)}, G^{(3)}, G^{(4)})$, but only for flexible ones can $M$ be a matching. The ultimate goal is to find and classify all matchings. We start from a given coupling $S_1$, and construct another coupling $S_2 = (g^{(3)}, g^{(4)}, h^{(3)}, h^{(4)})$ such that $Z(M)$ is an infinite set.

## 3 REDUCIBILITY AND SINGULARITY IN $\mathbb{C}P^1$

In Definition 2.6, the reason we consider the zero set in $\mathbb{C}P^1$ is that $x_i, y_j$ actually represent tangent values of angles (see equation (3), so it is possible for them to be evaluated as $\infty$. $g^{(1)} = 0$ implies

$$a_i y_i^2 + b_i + c_i y_i^2 \bar{x}_i^2 - \bar{x}_i y_i + e_i \bar{x}_i^2 = 0, \;\; \bar{x}_i = \tan\left(\frac{\alpha_i + \pi}{2}\right) = -\frac{1}{x_1}. \tag{7}$$

When $x_i = \infty$, i.e. $\bar{x}_i = 0$, the value(s) of $y_i$ can be recovered from $a_i y_i^2 + b_i = 0$. Conversely, when $a_i y_i^2 + b_i = 0$, $x_i = \frac{c_i y_i^2 + e_i}{-y_i}$ or $\infty$.

**Example 3.1.** Let us consider a polynomial

$$g^{(1)} = c_1 y_1^2 + x_1 y_1 + e_1, \;\; c_1 e_1 \neq 0,$$

which is not factorizable because, if so, $g^{(1)} = f(y_1)h(x_1, y_1)$ and, for constant $c \in \mathbb{C}$ such that $f(c) = 0$, $g^{(1)}(x_1, c) = f(c)h(x_1, c) \equiv 0 \in \mathbb{C}[x_1]$, which further implies $c = 0$ and $e_1 = 0$, a contradiction! It seems that in $Z(g^{(1)}) \subset (\mathbb{C}P^1)^2$, the value of $x_1$ is uniquely determined by $x_1 = \frac{c_1 y_1^2 + e_1}{-y_1}$. However, this is not true. Lemma 2.3 refers $g^{(1)}$ to an antideltoid in which $\lambda_1 + \gamma_1 = \mu_1 + \delta_1 = \pi$. If the layout of the antideltoid is as described in Fig. 4, then $y_1$ is free to change while $x_1$ remains $\infty$. One can also easily see this by rewriting $g^{(1)}$ into the form of equation (7).

[5]We are going to apply the notation $Z(\cdot)$ on different ideals, the dimension depends on the number of variables contained by the ideal, e.g. $Z(S) \subset (\mathbb{C}P^1)^5$, $Z(g^{(i)}) \subset (\mathbb{C}P^1)^2$, etc...

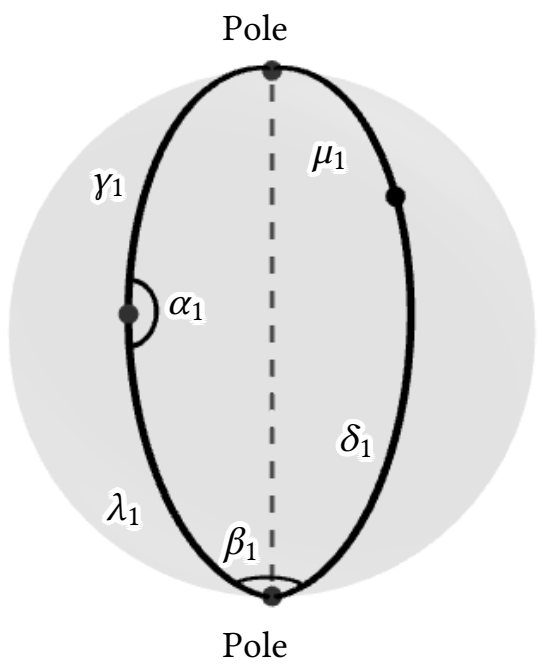

Fig. 4. When $(\lambda_1, \gamma_1)$ lie on the same great cycle and so do $(\mu_1, \delta_1)$, $\beta_1$ (i.e. $y_1$) is free to change while $\alpha_1$ remains $\pi$ (i.e. $x_1 = \tan(\frac{\pi}{2}) = \infty$).

The above example shows, even $g^{(1)}$ is irreducible, $Z(g^{(1)})$ is still reducible as a union of two curves $-x_1y_1 = c_1y_1^2 + e_1$ and $x_1 = \infty$. In fact, this inconsistency only occurs in antideltoids, and will be fixed in Section 4 where we convert all antideltoids to deltoids.

*Definition 3.1.* Let $P_i$ be the projection map of the $i$th coordinate in $(\mathbb{C}P^1)^n$ and $Z \subset (\mathbb{C}P^1)^n$. If $P_i^{-1}(c_0) \cap Z$ is an infinite set for some $c_0 \in \mathbb{C}P^1$, we call $P_i^{-1}(c_0) \cap Z$ a **constant branch** or a $c_0$-**branch** of $Z$.

*Definition 3.2.* $g^{(i)}$ (or $G^{(i)}$) is said to be **singular** if $a_ie_i = b_ic_i = 0$, it is **non-singular** otherwise. Similarly, a coupling (or matching) is said to be **singular** if it contains a singular $g^{(i)}$ (or $G^{(i)}$), it is **non-singular** otherwise.

Lemma 2.3 immediately shows every singular $g^{(i)}$ corresponds to an (anti)deltoid and vice versa.

Proposition 3.3. *$g^{(i)}$ is singular if and only if $Z(g^{(i)})$ has a constant branch.*

Proof. Recall equation (3) in which we have set $x_i = \tan\left(\frac{\alpha_i}{2}\right), y_i = \tan\left(\frac{\beta_i}{2}\right)$. Hence, $g^{(i)} = 0$ has four equivalent representations with substitutions

$$\bar{x}_i = \tan\left(\frac{\alpha_i + \pi}{2}\right) = -\frac{1}{x_1},\ \ \bar{y}_i = \tan\left(\frac{\beta_i + \pi}{2}\right) = -\frac{1}{y_1}.$$

$$g^{(i)} = 0 \Leftrightarrow \begin{cases} a_ix_i^2y_i^2 + b_ix_i^2 + c_iy_i^2 + x_iy_i + e_i = 0, \\ -a_iy_i^2 - b_i - c_iy_i^2\bar{x}_i^2 + \bar{x}_iy_i - e_i\bar{x}_i^2 = 0, \\ -a_ix_i^2 - b_ix_i^2\bar{y}_i^2 - c_i + x_i\bar{y}_i - e_i\bar{y}_i^2 = 0, \\ a_i + b_i\bar{y}_i^2 + c_i\bar{x}_i^2 + \bar{x}_i\bar{y}_i + e_i\bar{x}_i^2\bar{y}_i^2 = 0. \end{cases} \tag{8}$$

For conditions (iii)-(vi) of Lemma 2.3, it is easy to verify that $Z(g^{(i)})$ has either a 0-branch or an $\infty$-branch, so the singularity of $g^{(i)}$ guarantees $Z(g^{(i)})$ a constant branch. Conversely, for non-singular $g^{(i)}$, none of the conditions (iii)-(vi) hold. No matter which constant value we assign to $x_i$ or $y_i$ (not both), the above four representations can never be identically 0, which means $Z(g^{(i)})$ has no constant branch. □

### 3.1 Non-singular properties

Lemma 3.4. *A non-singular polynomial $g^{(i)}$ is reducible if and only if $a_i = e_i = 0$ or $b_i = c_i = 0$. Moreover, when $g^{(i)}$ is reducible, the factorization must be in the form*

$$g^{(i)} = \begin{cases} c_i(kx_i - y_i)(k'x_i - y_i) \text{ if } a_i = e_i = 0, \\ a_i(x_iy_i - k)(x_iy_i - k') \text{ if } b_i = c_i = 0, \end{cases}$$

*where $\{k, k'\} \cap \{0, \pm 1\} = \emptyset$.*

Proof. Firstly, Note that when $a_i = e_i = 0$ or $b_i = c_i = 0$, $g^{(i)}$ can be factorized as stated and $k, k' \neq 0$ is directly from the non-singularity of $g^{(i)}$. Also, $k, k' \neq \pm 1$, otherwise the second inequality of (6) fails. Now, if neither $a_i = e_i = 0$ nor $b_i = c_i = 0$, $g^{(i)}$ must be irreducible. We prove this by contradiction. Assume $g^{(i)}$ is reducible and none of $a_i = e_i = 0$ or $b_i = c_i = 0$ holds. Proposition 3.3 implies $g^{(i)}$ has no factors in $\mathbb{C}[x_i]$ or $\mathbb{C}[y_i]$. So any irreducible factor of $g^{(i)}$ must be in the form

$$px_iy_i + qx_i - ry_i - s, p \neq 0 \text{ or } q \neq 0.$$

This means, by regarding $x_i$ as the unknown and $y_i$ as a parameter, $g^{(i)} = 0$ admits a rational solution $x_i = \frac{ry_i+s}{py_i+q} \in \mathbb{C}(y_i)$. In other words, the discriminant

$$-4a_ic_iy_i^4 + (1 - 4a_ie_i - 4b_ic_i)y_i^2 - 4b_ie_i$$

is a square in $\mathbb{C}[y_i]$. With this information, we can claim that $a_ic_i \neq 0$. Otherwise, we must have $(1 - 4a_ie_i - 4b_ic_i) = 0$ or $b_ie_i = 0$. The first one fails given inequalities (6), and the second is also impossible due to non-singular $g^{(i)}$. So the discriminant can be factorized into $-4a_ic_i(y_i^2 - l)(y_i^2 - l')$ where $l = l'$ since it is a square. Consequently,

$$-4a_ic_iy_i^2 + (1 - 4a_ie_i - 4b_ic_i)y_i - 4b_ie_i = -4a_ic_i(y_i - l)^2$$

is a square in $\mathbb{C}[y_i]$ with a zero discriminant $(1 - 4a_ie_i - 4b_ic_i)^2 - 64a_ib_ic_ie_i = 0$, contrary to inequalities (6). □

Combine Lemma 2.3 and 3.4, one can conclude the following without a proof.

Theorem 3.5. *Given an ideal $M = (G^{(1)}, G^{(2)}, G^{(3)}, G^{(4)}) \subset \mathbb{C}[x_1, x_2, x_3, x_4]$ such that $Z(M)$ has no constant branch, then $M$ is a matching if and only if*

$$\gcd(Res(G^{(1)}, G^{(2)}; x_2), Res(G^{(3)}, G^{(4)}; x_4)) \neq 1.$$

Proof. Set $R_1 = \text{Res}(G^{(1)}, G^{(2)}; x_2), R_3 = \text{Res}(G^{(3)}, G^{(4)}; x_4)$. Given $Z(M)$ has no constant branch, as a subset of $(\mathbb{C}P^1)^4$, $Z(M)$ is an infinite set if and only if $Z(M) \cap \mathbb{C}^4$ is an infinite set. On the other hand, the $(x_1, x_3)$-projection of $Z(G^{(1)}, G^{(2)})$ and $Z(G^{(3)}, G^{(4)})$ lie in $Z(R_1)$ and $Z(R_3)$ respectively. Thus two algebraic curves $R_1 = 0, R_3 = 0$ have overlapping branches in $\mathbb{C}^2$, i.e. $\gcd(R_1, R_3) \neq 1$. □

## 4 BASIC REDUCTIONS

The purpose of this section is to considerably reduce our future discussion.

*Definition 4.1.* Let $P_4$ be the permutation group generated by $(1234)$, couplings

$$S_1 = (g^{(1)}, h^{(1)}, g^{(2)}, h^{(2)}),\ S_2 = (g^{(3)}, h^{(3)}, g^{(4)}, h^{(4)});$$

$$\tilde{S}_1 = (\tilde{g}^{(1)}, \tilde{h}^{(1)}, \tilde{g}^{(2)}, \tilde{h}^{(2)}),\ \tilde{S}_2 = (\tilde{g}^{(3)}, \tilde{h}^{(3)}, \tilde{g}^{(4)}, \tilde{h}^{(4)}).$$

$(\tilde{S}_1, \tilde{S}_2)$ is called a **relabeling** of $(S_1, S_2)$ if $\exists \sigma \in P_4$ such that, $\forall i \in \{1,2,3,4\}$, $\tilde{g}^{(i)} \equiv g^{(\sigma(i))}$ and $\tilde{h}^{(i)} \equiv h^{(\sigma(i))}$ under substitution $(x_{\sigma(i)}, y_{\sigma(i)}) = (x_i, y_i)$. In other words, $\tilde{g}^{(i)}$ (or $\tilde{h}^{(i)}$) is obtained from $g^{(\sigma(i))}$ (or $h^{(\sigma(i))}$) by replacing the variables $(x_{\sigma(i)}, y_{\sigma(i)})$ by $(x_i, y_i)$ but keeping the corresponding coefficients $(a_{\sigma(i)}, b_{\sigma(i)}, ...)$ etc.

It is obvious that relabeling will not affect the cardinality of the fiber product given

$$Z(S_1) \underset{\{x_1,x_3\}}{\times} Z(S_2) \cong Z(g^{(2)}, h^{(2)}, g^{(3)}, h^{(3)}) \underset{\{x_2,x_4\}}{\times} Z(g^{(4)}, h^{(4)}, g^{(1)}, h^{(1)}).$$

Relabeling is just rearranging of polynomials so that $Z(S_i)$ would be easier to analyze.

Another reduction is on antiisograms and antideltoids, whose construction can be converted to their counterparts. Suppose we have a mesh with an antiisogram. With a relabeling, in its corresponding matching $M$ we have $a_1 = e_1 = 0$ (Lemma 2.3). Set

$$M' := (g^{(1)}, h^{(1)}, g^{(2)}, h^{(2)}, g^{(3)}, h^{(3)}, g^{(4)}, h^{(4)}). \tag{9}$$

Notice that $Z(M) \cong Z(M')$ due to the birational map

$$\begin{aligned} f:\ & Z(M) \leftrightarrow Z(M'); \\ & (x_1, ..., x_4) \leftrightarrow (x_1, y_1, ..., x_4, y_4) = \left(x_1, \tfrac{(1-f_1^2)x_2+2f_1}{1-f_1^2-2f_1x_2}, ..., x_4, \tfrac{(1-f_4^2)x_1+2f_4}{1-f_4^2-2f_4x_1}\right). \end{aligned} \tag{10}$$

Set

$$\bar{g}^{(1)}(x_1, \bar{y}_1) = \bar{a}_1 x_1^2 \bar{y}_1^2 + \bar{e}_1 + x_1\bar{y}_1,\ \bar{h}^{(1)}(\bar{y}_1, x_2) = (1 - \bar{f}_1^2)(\bar{y}_1 - x_2) - 2\bar{f}_1(\bar{y}_1 x_2 + 1)$$

where $(\bar{a}_1, \bar{e}_1) = -(b_1, c_1)$ and $\bar{f}_1 = \frac{f_1+1}{1-f_1}$ if $f_1 \le 0$, $\bar{f}_1 = \frac{f_1-1}{1+f_1}$ if $f_1 > 0$. Then set

$$\bar{M} := (\bar{g}^{(1)}, \bar{h}^{(1)}, g^{(2)}, h^{(2)}, g^{(3)}, h^{(3)}, g^{(4)}, h^{(4)}).$$

The substitutions above equation (8) induce a birational map between $Z(M')$ and $Z(\bar{M})$ thus $Z(M) \cong Z(\bar{M})$. And in $\bar{M}$ we have $\bar{g}^{(1)}$ correspond to an isogram. If there are more antiisograms, we can continue to work out more $\bar{g}^{(i)}$. This trick is obviously able to convert antideltoids to deltoids as well.

## 5 NON-SINGULAR FLEXIBLE MESHES WITH ISOGRAMS

The meshes to be dealt with in this section correspond to **isogonal matchings** that consist of non-singular but reducible $g^{(i)}$ (Lemma 2.3 and 3.4). This particular type of mesh is the first non-trivial flexible mesh with skew faces that has ever been constructed [6]. The author has built several flexible meshes of 'isogonal type', and before him, the existence of a flexible mesh with skew faces was unknown. Hereby, we will complete his work by showing a full construction of all possible isogonal matchings.

Theorem 5.1. *Given a matching $M = (G^{(1)}, G^{(2)}, G^{(3)}, G^{(4)})$ in which $a_i e_i \neq 0$ and $b_i = c_i = 0$ for $i = 1, 2, 3, 4$. Set*

$$k_i = \frac{-1 \pm \sqrt{1 - 4a_i e_i}}{2a_i},\ N_i = \begin{pmatrix} -2f_i & k_i(1 - f_i^2) \\ 1 - f_i^2 & 2k_i f_i \end{pmatrix}.$$

*Consider all combinations of signs '±' that appear in $k_i$, the mesh is flexible if and only if there exists a combination such that $N_4 N_3 N_2 N_1$ is a scalar matrix.*

We refer to [6] for visualizations of such flexible meshes.

Proof. By Lemma 3.4, the corresponding $g^{(i)}$ can be factorized as

$$g^{(i)} = a_i x_i^2 y_i^2 + x_i y_i + e_i = a_i(x_i y_i - k_i)(x_i y_i - k_i')$$

where $\{k_i, k_i'\} = \left\{\frac{-1\pm\sqrt{1-4a_ie_i}}{2a_i}\right\}$. Consider the coupling $S_1 = (g^{(1)}, h^{(1)}, g^{(2)}, h^{(2)})$, $Z(S_1)$ has four components of the same format, e.g.

$$W_1 := Z(x_1y_1 - k_1, h^{(1)}, x_2y_2 - k_2, h^{(2)}).$$

It is clear that $W_1$ is a curve that can be iteratively parametrized by $x_1 \in \mathbb{C}P^1$,

$$\begin{aligned}\omega(x_1) &= (x_1, y_1(x_1), x_2(y_1), y_2(x_2), x_3(y_2)) \\ &= \left(x_1, \frac{k_1}{x_1}, \frac{(1-f_1^2)y_1 - 2f_1}{2f_1y_1 + 1 - f_1^2}, \frac{k_2}{x_2}, \frac{(1-f_2^2)y_2 - 2f_2}{2f_2y_2 + 1 - f_2^2}\right).\end{aligned}$$

If we formally express the above relations in a matrix form

$$N = \begin{pmatrix} u & v \\ s & t \end{pmatrix} : \mathbb{C}P^1 \to \mathbb{C}P^1;\ x \mapsto N(x) := \frac{ux+v}{sx+t},\ ut - sv \neq 0,$$

this is the so-called Möbius transformation on $\mathbb{C}P^1$ which is compatible with matrix multiplications. e.g.

$$x_{i+1} = \begin{pmatrix} 1-f_i^2 & -2f_i \\ 2f_i & 1-f_i^2 \end{pmatrix}(y_i) = \begin{pmatrix} 1-f_i^2 & -2f_i \\ 2f_i & 1-f_i^2 \end{pmatrix}\begin{pmatrix} 0 & k_i \\ 1 & 0 \end{pmatrix}(x_i) = N_i(x_i).$$

So we have $x_3 = N_2N_1(x_1)$. Likewise, in a component $W_2$ of coupling $S_2$ we have $x_1 = N_4N_3(x_3)$. Since $Z(M) \cong Z(S_1) \times_{\{x_1,x_3\}} Z(S_2)$ is an infinite set, there must exist components $W_1, W_2$ of $S_1, S_2$, respectively, such that $W_1 \times_{\{x_1,x_3\}} W_2$ is an infinite set. Therefore, $\{x_3 = N_2N_1(x_1), x_1 = N_4N_3(x_3)\}$ have infinitely many solutions, in other words $N_4N_3N_2N_1(x_1) = x_1$. □

Considering Bricard equation (4), $k_i$ in Theorem 5.1 can be also expressed as $\frac{\sin(\gamma_i-\mu_i)}{\sin(\gamma_i)\pm\sin(\mu_i)}$ which is more computational friendly to mesh construction. However, according to Remark 1, once we have the coefficients $(a_i, b_i, c_i, e_i, f_i)$, the process of generating a mesh is elementary. Thus, we now only focus on determining coefficients, or equivalently, matchings.

**Example 5.1** (General case for isogonal matching). Set

$$(a_1, e_1, a_2, e_2, a_3, e_3, a_4, f_1, f_2) \tag{11}$$

as free parameters such that all $a_i, e_i$ nonzero and all $f_i \in (-1, 1]$. $(e_4, f_3, f_4)$ are to be determined. For $i = 1, 2, 3$ set $k_i = \frac{-1\pm\sqrt{1-4a_ie_i}}{2a_i}$, and for $i = 1, 2$ set

$$N_i = \begin{pmatrix} -2f_i & k_i(1-f_i^2) \\ 1-f_i^2 & 2k_if_i \end{pmatrix}. \tag{12}$$

Let

$$A = \frac{2(us+vt)}{s^2+t^2-u^2-v^2} \text{ where } \begin{pmatrix} u & v \\ s & t \end{pmatrix} = \begin{pmatrix} 0 & k_3 \\ 1 & 0 \end{pmatrix} N_2N_1, \tag{13}$$

and let $F_3 \in \mathbb{R}P^1, f_3 \in (-1, 1]$ be such that $\frac{2F_3}{1-F_3^2} = A$, $\frac{2f_3}{1-f_3^2} = F_3$. Then set

$$\begin{pmatrix} u' & v' \\ s' & t' \end{pmatrix} = \begin{pmatrix} 1-f_3^2 & -2f_3 \\ 2f_3 & 1-f_3^2 \end{pmatrix}\begin{pmatrix} u & v \\ s & t \end{pmatrix}, \tag{14}$$

and define $(k_4, F_4)$ to be $(\frac{v'}{s'}, \frac{t'}{s'})$ if $s' \neq 0$, or $(-\frac{u'}{t'}, \infty)$ otherwise, and set $f_4 \in (-1, 1]$ be such that $\frac{2f_4}{1-f_4^2} = F_4$. Adjust the value of $a_4$ so that $e_4 = -a_4 k_4^2 - k_4 \neq 0$. All of the above signs '±' are free to choose. Start over if inequalities (6) do not hold for all $i$.

Proposition 5.2. *Example 5.1 contains all isogonal matchings.*

Proof. Please note the setups for $k_i$ and $N_i$ in Example 5.1 coincide with those in Theorem 5.1. So it is sufficient to prove that Example 5.1 provides all solutions for $N_4N_3N_2N_1$ being scalar. We start from arbitrary coefficients (11) and show that the remaining $(e_4, f_3, f_4)$ must satisfy their determinations if $N_4N_3N_2N_1$ is scalar. For Möbius transformations, any scalar matrix is the identity map

$$\begin{pmatrix} n & 0 \\ 0 & n \end{pmatrix}(x) = \frac{nx + 0}{0x + n} = x, \forall n \neq 0.$$

Hence we may regard $N_i \in \mathrm{PGL}(2, \mathbb{R})$ and the inverse of a matrix, despite a scalar, is just its adjoint matrix. Simple observation shows

$$N_i \in T := \left\{ \begin{pmatrix} u & v \\ s & t \end{pmatrix} : uv = -st \right\},\ N_i^{-1} \in T' := \left\{ \begin{pmatrix} u' & v' \\ s' & t' \end{pmatrix} : u's' = -v't' \right\},$$

which implies $N_3N_2N_1 \in T'$ since $N_4N_3N_2N_1$ is supposed to be scalar. Notice that

$$f:\ (-1, 1] \to \mathbb{R}P^1;\ x \mapsto \frac{2x}{1 - x^2}$$

is onto given $f(\tan(\alpha)) = \tan(2\alpha)$. Hence, the $f_3$ below equation (13) does exist. This is equivalent to saying that $\exists f_3 \in (-1, 1]$ such that matrix (14) belongs to $T'$, i.e. $u's' + v't' = 0$, and we must have $s' \neq 0$ or $t' \neq 0$ due to $N_i$ invertible. Thus, in $\mathrm{PGL}(2, \mathbb{R})$, we have

$$N_4 = \begin{pmatrix} t' & -v' \\ -s' & u' \end{pmatrix} = \begin{pmatrix} -\frac{t'}{s'} & \frac{v'}{s'} \\ 1 & -\frac{u'}{s'} \end{pmatrix} \text{ or } N_4 = \begin{pmatrix} -1 & 0 \\ 0 & -\frac{u'}{t'} \end{pmatrix}.$$

The reason of normalizing $N_4$ into the above form is to recover the value of $k_4$ and $F_4 = \frac{2f_4}{1-f_4^2}$ from the definition of $N_4$ in Theorem 5.1:

$$N_4 = \begin{pmatrix} -2f_4 & k_4(1 - f_4^2) \\ 1 - f_4^2 & 2k_4 f_4 \end{pmatrix} = \begin{cases} \begin{pmatrix} -F_4 & k_4 \\ 1 & k_4 F_4 \end{pmatrix} & \text{if } f_4 \in (-1, 1), \\ \begin{pmatrix} -1 & 0 \\ 0 & k_4 \end{pmatrix} & \text{if } f_4 = 1. \end{cases}$$

Finally, $e_4 = -a_4 k_4^2 - k_4 \neq 0$ given how $k_i$ was defined in Theorem 5.1. □

## 6 SINGULAR FLEXIBLE MESHES WITH CONSTANT BRANCH

For singular matching, we first consider the easiest case, **constant matchings**, which means $Z(M)$ has a constant branch (Definition 3.1). It is unnecessary to consider antideltoids due to Section 4.

Proposition 6.1. *For every constant matching $(G^{(1)}, G^{(2)}, G^{(3)}, G^{(4)})$, there are $G^{(i)}$ and $G^{(j)}$, $i \neq j$, such that $b_i = e_i = 0$ and $c_j = e_j = 0$.*

Proof. From equations (9) and (10) we have $Z(M) \cong Z(M')$ so that $Z(M')$ admits a constant branch, say $W$. By symmetry, we may assume $x_1$ varies and $x_2$ is constant in $W$. Hence $y_1 = \frac{(1-f_1^2)x_2+2f_1}{1-f_1^2-2f_1x_2}$ is also constant. This means $Z(g^{(1)})$ has a constant branch. Since antideltoids are excluded, by Proposition 3.3, $g^{(1)}$ must satisfy condition (v) of Lemma 2.3, i.e. $b_1 = e_1 = 0$, so that $y_1 \equiv 0$. Now, in $W$, we have $x_1$ varies and $x_2$ fixed, so there exists $j \in \{2,3,4\}$ such that $x_j$ is fixed and $x_{j+1}$ varies ($x_5 := x_1$). The same reasoning shows $c_j = e_j = 0$. □

Without loss of generality, in a constant branch, we may always assume that $x_1$ varies but $x_2$ is fixed, the above proposition actually reveals a systematic construction.

**Example 6.1** (General case for constant matching)**.** Choose $j \in \{2,3,4\}$, set $(b_1, e_1, c_j, e_j)$ all 0 and $\{a_1, c_1, f_4, a_j, b_j\}$ arbitrary. For $i \notin \{1, j\}$, arbitrarily set $\{a_i, b_i, c_i, e_i, f_{i-1}\}$. Set $Y_1 = 0$, and for $1 < k < j$, if $k$ exists, iteratively set

$$X_k := \frac{(1-f_{k-1}^2)Y_{k-1} - 2f_{k-1}}{1 - f_{k-1}^2 + 2f_{k-1}Y_{k-1}}, \quad Y_k := \frac{-1 \pm \sqrt{1 - 4(a_kX_k^2 + c_k)(b_kX_k^2 + e_k)}}{2(a_kX_k^2 + c_k)}.^6$$

Start over if inequalities (6) do not hold for all $i$ or $Y_k$ is not real. Finally, set $f_{j-1} \in (-1, 1]$ such that $\frac{2f_{j-1}}{1-f_{j-1}^2} = Y_{j-1}$. It is easy to check that this matching has a constant branch in which $x_k \equiv X_k$ for all $2 \le k \le j$. In particular, $x_j \equiv 0$.

## 7 SINGULAR FLEXIBLE MESHES WITH NON-CONSTANT BRANCH

Even though we do not need to consider antideltoids, a singular non-constant matching may contain one to four singular $G^{(i)}$ satisfying condition (iii) or (v) of Lemma 2.3 (the remaining $G^{(i)}$ satisfy (i)), and how these heterogeneous $G^{(i)}$ are arranged in order also matters. The discussion for such matchings seems overwhelming. However, most of the combinations of $G^{(i)}$ cannot form a matching ($Z(M)$ is a finite set), and remarkably, there are essentially two main cases only!

### 7.1 Preparations

Recall equations (9) and (10), for a matching $M$ in which $g^{(1)} = a_1x_1^2y_1^2 + b_1x_1^2 + x_1y_1 = x_1(a_1x_1y_1^2 + b_1x_1 + y_1)$, we have

$$Z(M) \cong Z(g^{(1)}, h^{(1)}, \ldots) = Z(x_1, h^{(1)}, \ldots) \cup Z(a_1x_1y_1^2 + b_1x_1 + y_1, h^{(1)}, \ldots). \tag{15}$$

Since constant branches are now out of interest, we only need to investigate the last component. In addition, we may assume $a_1b_1 \neq 0$. The reasons are as follows.

- If $(a_1, b_1)$ are both 0, $Z(M) \cong Z(x_1, \ldots) \cup Z(y_1, \ldots)$ which only contains constant branches, a contradiction;
- If one of $(a_1, b_1)$ is 0, the non-constant branch of $Z(M)$ is in the form $Z(kx_1 - y_1, \ldots)$ or $Z(x_1y_1 - k, \ldots)$ in which $g^{(1)}$ contributes as a non-singular polynomial (Lemma 3.4).

In this regard, we introduce some customized notations for this section.

*Definition 7.1.* All irreducible $g^{(i)}$ are denoted by $g_0^{(i)}$. $g_1^{(i)}, g_3^{(i)}, g_5^{(i)} \in \mathbb{R}[x_i, y_i]$ are polynomials in the form

$$g_1^{(i)} := x_iy_i - k_i, \; g_3^{(i)} := a_ix_iy_i^2 + b_ix_i + y_i, \; g_5^{(i)} := a_ix_i^2y_i + c_iy_i + x_i$$

where the coefficients are nonzero and satisfy

$$k_i \notin \{0, \pm 1\}, \; |b_i - a_i| < 1, \; |c_i - a_i| < 1. \tag{16}$$

For $j \in \{0, 1, 3, 5\}$, $G_j^{(i)} \in \mathbb{R}[x_i, x_{i+1}]$ is the resultant

$$G_j^{(i)}(x_i, x_{i+1}) := \mathrm{Res}(g_j^{(i)}, h^{(i)}; y_i).$$

[6] More precisely, $y_k = Y_k$ is a root of $g^{(k)}(X_k, y_k) = 0$ ($a_k = c_k = 0$ could happen).

The index $i$ is always considered modulo 4, e.g., $x_5 = x_1$.

The reader should realize that, for $j \in \{1,3,5\}$, $g_j^{(i)}$ (or $G_j^{(i)}$) defined above are prime factors of $g^{(i)}$ (or $G^{(i)}$) when condition (i), (iii), or (v) of Lemma 2.3 are respectively satisfied, and the inequalities are inherited from Lemma 3.4 and Definition 2.4.

*Definition 7.2.* Let $S = (g^{(i)}, h^{(i)}, g^{(i+1)}, h^{(i+1)})$ be a coupling. A reduced coupling $S'$ (of $S$) is in the form

$$S' = (g_{j_i}^{(i)}, h^{(i)}, g_{j_{i+1}}^{(i+1)}, h^{(i+1)})$$

where $j_i, j_{i+1} \in \{0,1,3,5\}$, $g_{j_i}^{(i)}, g_{j_{i+1}}^{(i+1)}$ are prime factors of $g^{(i)}, g^{(i+1)}$ respectively. Further, set $r_{S'} \in \mathbb{C}[x_i, y_{i+1}]$ and $R_{S'} \in \mathbb{C}[x_i, x_{i+2}]$ be the resultants

$$r_{S'} := \mathrm{Res}(G_{j_i}^{(i)}, g_{j_{i+1}}^{(i+1)}; x_{i+1}),\ \ R_{S'} := \mathrm{Res}(G_{j_i}^{(i)}, G_{j_{i+1}}^{(i+1)}; x_{i+1}).$$

The index $i$ is always considered modulo 4, e.g., $x_5 = x_1$.

It is easy to see that $Z(S')$ is simply $Z(S)$ kicking out its constant branches. So, if a matching has a non-constant branch, we only need two reduced couplings $S_1', S_2'$ such that $Z(S_1') \times_{\{x_1,x_3\}} Z(S_2')$ is an infinite set. This is equivalent to saying that $\gcd(R_{S_1'}, R_{S_2'}) \neq 1$ according to Theorem 3.5. To take advantage of this, we need to develop some knowledge to analyze the factorization of $R_{S_i'}$.

*Definition 7.3.* For $p(x,y) \in \mathbb{C}[x,y]$, let $\deg_x(p)$ and $\deg_y(p)$ respectively be the highest degree of $x$ and $y$ in $p$. The **double degree** of $p(x,y)$ is

$$\deg_y^x(p) := (\deg_x(p), \deg_y(p)).$$

*Definition 7.4.* $\mathcal{M} : \mathbb{C}[x,y] \to \mathbb{C}[w,z]$ is called a **Möbius map** if there are matrices $\begin{pmatrix} u & v \\ s & t \end{pmatrix}, \begin{pmatrix} U & V \\ S & T \end{pmatrix} \in \mathrm{GL}_2(\mathbb{C})$ such that for $(m,n) = \deg_y^x(p)$,

$$\mathcal{M}(p(x,y)) = (sw+t)^m (Sz+T)^n p\left(\frac{uw+v}{sw+t}, \frac{Uz+V}{Sz+T}\right).$$

Lemma 7.5. *Let $\mathcal{M} : \mathbb{C}[x,y] \to \mathbb{C}[w,z]$ be a Möbius map. For all $p \in \mathbb{C}[x,y]$ such that $p$ has no factors in $\mathbb{C}[x]$ or $\mathbb{C}[y]$, $\mathcal{M}(p)$ has no factors in $\mathbb{C}[w]$ or $\mathbb{C}[z]$ and*

$$deg_z^w(\mathcal{M}(p)) = deg_y^x(p).$$

Proof. $\forall f(x,y) \in \mathbb{C}[x,y], \forall g(w,y) \in \mathbb{C}[w,y]$, define

$$\mathcal{M}_1(f) := (sw+t)^{\deg_x(f)} f\left(\frac{uw+v}{sw+t}, y\right),\ \ \mathcal{M}_2(g) = (Sz+T)^{\deg_y(g)} g\left(w, \frac{Uz+V}{Sz+T}\right).$$

Clearly, we have $\mathcal{M} = \mathcal{M}_2 \circ \mathcal{M}_1$. Let $P = \mathcal{M}_1(p)$, $m = \deg_x(p)$ and $m' = \deg_w(P)$. It is sufficient to prove $m = m'$ and $P$ has no factors in $\mathbb{C}[w]$. $m' \leq m$ is obvious. Conversely, direct calculation shows $(u-sx)^m P\left(\frac{tx-v}{u-sx}, y\right) = (ut-sv)^m p(x,y)$. Since $(u-sx)$ cannot be a factor of $p$, $m' = m$ holds. Now, if $P$ has a factor $f(w)$, then $P = f(w) \cdot P'$ where $\deg_w(P') < m$. Define Möbius map

$$\mathcal{M}_1'(g(w,y)) = (u-sx)^{\deg_w(g)} g\left(\frac{tx-v}{u-sx}, y\right).$$

We have $\mathcal{M}_1'(P) = \mathcal{M}_1'(f) \cdot \mathcal{M}_1'(P') = (ut-sv)^m p(x,y)$ where $\mathcal{M}_1'(f) \in \mathbb{C}[x]$ and $\deg_x(\mathcal{M}_1'(P')) < m$. Given $p$ has no factors in $\mathbb{C}[x]$, $\mathcal{M}_1'(f)$ must be a constant so that $\deg_x(\mathcal{M}_1'(P')) = \deg_x(p) = m$, a contradiction! □

| $j$ | 0 | 1 | 3 | 5 |
|---|---|---|---|---|
| $deg^{x_i}_{x_{i+1}}(G^{(i)}_j)$ | $(2,2)$ | $(1,1)$ | $(1,2)$ | $(2,1)$ |

Lemma 7.6. *In Definition 7.1 we have* $deg^{x_i}_{y_i}(g^{(i)}_j) = deg^{x_i}_{x_{i+1}}(G^{(i)}_j)$. *In particular,*

Proof. Set Möbius map

$$\mathcal{M}_i(g(x,y)) = (1 - f_i^2 - 2f_i x_{i+1})^{\deg_y(g)} g\left(x, \frac{(1-f_i^2)x_{i+1} + 2f_i}{1 - f_i^2 - 2f_i x_{i+1}}\right). \tag{17}$$

Replace $(x, y)$ by $(x_i, y_i)$ one can easily verify that $G^{(i)}_j = \mathcal{M}_i(g^{(i)}_j)$. □

Theorem 7.7. *Let* $S' = (g^{(1)}_{j_1}, h^{(1)}, g^{(2)}_{j_2}, h^{(2)})$ *be a reduced coupling in which* $deg^{x_1}_{x_2}(g^{(1)}_{j_1}) = (m_1, n_1)$, $deg^{x_2}_{x_3}(g^{(2)}_{j_2}) = (m_2, n_2)$. *Then*

$$deg^{x_1}_{y_2}(r_{S'}) = deg^{x_1}_{x_3}(R_{S'}) = (m_1 m_2, n_1 n_2).$$

*Moreover,* $r_{S'}$ *is reducible if and only if* $R_{S'}$ *is reducible.*

Proof. First, notice that both $r_{S'}$ and $R_{S'}$ have no factors in $\mathbb{C}[x_1], \mathbb{C}[y_2]$ or $\mathbb{C}[x_3]$ because, otherwise, $Z(S')$ would have a constant branch. Let $(m, n) = \deg^{x_1}_{y_2}(r_{S'})$ and set $i = 2$ in Möbius map (17), we have

$$\begin{aligned}
\mathcal{M}_2(r_{S'}(x_1, y_2)) &= (1 - f_2^2 - 2f_2x_3)^n r_{S'}\left(x_1, \frac{(1-f_2^2)x_3 + 2f_2}{1 - f_2^2 - 2f_2x_3}\right) \\
&= (1 - f_2^2 - 2f_2x_3)^n \mathrm{Res}\left(G^{(1)}_{j_1}, g^{(2)}_{j_2}\left(x_2, \frac{(1-f_2^2)x_3 + 2f_2}{1 - f_2^2 - 2f_2x_3}\right); x_2\right) \\
\text{Lemma 7.6 proof} \rightarrow &= (1 - f_2^2 - 2f_2x_3)^n \mathrm{Res}\left(G^{(1)}_{j_1}, \frac{G^{(2)}_{j_2}}{(1 - f_2^2 - 2f_2x_3)^{n_2}}; x_2\right) \\
&= (1 - f_2^2 - 2f_2x_3)^{n - n_1 n_2} \mathrm{Res}\left(G^{(1)}_{j_1}, G^{(2)}_{j_2}; x_2\right) \\
&= (1 - f_2^2 - 2f_2x_3)^{n - n_1 n_2} R_{S'} \in \mathbb{C}[x_1, x_3].
\end{aligned}$$

By Lemma 7.5, $\mathcal{M}_2(r_{S'})$ has no factors in $\mathbb{C}[x_3]$, so it must be $n = n_1 n_2$. Consequently, $\deg^{x_1}_{y_2}(R_{S'}) = \deg^{x_1}_{y_2}(r_{S'}) = (m, n_1 n_2)$. To prove $m = m_1 m_2$, we need an 'auxiliary coupling' $\tilde{S'} = (h^{(0)}, g^{(1)}_{j_1}, h^{(1)}, g^{(2)}_{j_2})$ and for $i = 1, 2$ define

$$\tilde{G}^{(i)}_j = \mathrm{Res}(g^{(i)}_{j_i}, h^{(i-1)}; x_i), \;\; \tilde{r}_{\tilde{S'}} = \mathrm{Res}(g^{(1)}, \tilde{G}^{(2)}; y_1), \;\; \tilde{R}_{\tilde{S'}} = \mathrm{Res}(\tilde{G}^{(1)}, \tilde{G}^{(2)}; y_1).$$

Now $x_1$ plays the same role in $(\tilde{r}_{\tilde{S'}}, \tilde{R}_{\tilde{S'}})$ as $y_2$ did in $(r_{S'}, R_{S'})$, and it is obvious that $\deg_{x_1}(r_{S'}) = \deg_{x_1}(\tilde{r}_{\tilde{S'}})$. Thus, a similar argument can show $m = m_1 m_2$.

Finally, it is clear that the factorizations of $r_{S'}$ and $R_{S'}$ induce each other through Möbius maps, hence they share the same reducibility. □

Apply Theorem 7.7 on Lemma 7.6, we immediately have the following result.

Corollary 7.8. *For reduced couplings* $S_1' = (g_{j_1}^{(1)}, h^{(1)}, g_{j_2}^{(2)}, h^{(2)})$ *and* $S_2' = (g_{j_3}^{(3)}, h^{(3)}, g_{j_4}^{(4)}, h^{(4)})$, $deg_{x_3}^{x_1}(R_{S_1'})$ *and* $deg_{x_3}^{x_1}(R_{S_2'})$ *are as follows.*

| $\{j_1, j_2\}/\{j_3, j_4\}$ | $\{1,1\}$ | $\{1,3\}$ | $\{1,5\}$ | $\{3,3\}$ | $\{5,5\}$ | $\{3,5\}$ |
|---|---|---|---|---|---|---|
| $deg_{x_3}^{x_1}(R_{S_1'})$ | $(1,1)$ | $(1,2)$ | $(2,1)$ | $(1,4)$ | $(4,1)$ | $(2,2)$ |
| $deg_{x_3}^{x_1}(R_{S_2'})$ | $(1,1)$ | $(2,1)$ | $(1,2)$ | $(4,1)$ | $(1,4)$ | $(2,2)$ |

Theorem 7.9. *Let* $S_1' = (g_{j_1}^{(1)}, h^{(1)}, g_{j_2}^{(2)}, h^{(2)})$ *and* $S_2' = (g_{j_3}^{(3)}, h^{(3)}, g_{j_4}^{(4)}, h^{(4)})$ *be reduced couplings such that* $j_i \in \{1,3,5\}$. *If* $Z(S_1') \times_{\{x_1,x_3\}} Z(S_2')$ *is an infinite set, up to a relabeling,* $\{j_1, j_2\}$ *and* $\{j_3, j_4\}$ *must be one of the following combinations.*

| | | |
|---|---|---|
| $\{j_1, j_2\}$ | $\{1,3\}$ | $\{3,5\}$ |
| $\{j_3, j_4\}$ | $\{1,5\}$ | $\{3,5\}$ |

Proof. According to Theorem 3.5 we have $\gcd(R_{S_1'}, R_{S_2'}) \neq 1$. In particular, when $\deg_{x_3}^{x_1}(R_{S_1'}) = \deg_{x_3}^{x_1}(R_{S_2'})$, $R_{S_1'}$ and $R_{S_2'}$ could be identical up to a constant. After a case-by-case consideration of Corollary 7.9, and relabel the polynomials if necessary, we obtain the above table. As for $\deg_{x_3}^{x_1}(R_{S_1'}) \neq \deg_{x_3}^{x_1}(R_{S_2'})$, remember that $R_{S_i'}$ has no factors in $\mathbb{C}[x_1]$ or $\mathbb{C}[x_3]$. Hence, for instance, $R_{S_1'}$ has no factor of double degree $(1,2)$ if $\deg_{x_3}^{x_1}(R_{S_1'}) = (2,2)$. So, the only possible combination left is $\{j_1, l_1\} = \{1,1\}$ and $\{j_2, l_2\} = \{3,5\}$, which goes back to $\{1,3\}$ and $\{1,5\}$ after relabeling. □

The above theorem reveals all types of reducible singular matchings with a non-constant branch. And according to Theorem 7.7, it is equivalent (and easier) to consider the reducibility of $r_{S'}$ rather than $R_{S'}$, so the following two lemmas arise.

Lemma 7.10. *For a reduced coupling* $S' = (g_3^{(1)}, h^{(1)}, g_5^{(2)}, h^{(2)})$, $r_{S'}$ *is reducible if and only if one of the following systems holds.*

- $\{b_1 = -a_1, c_2 = -a_2\}$; • $\{a_1c_2 = a_2b_1, f_1 = 0\}$; • $\{a_1a_2 = b_1c_2, f_1 = 1\}$.

*In particular, the factorizations of* $r_{S'}$ *are respectively as follows.*

- $r_{S'} = -[2f_1(f_1^2 - 1)(4a_1a_2x_1y_2 + 1) + ((f_1^2 - 1)^2 - 4f_1^2)(a_1x_1 - a_2y_2)]^2$;
- $r_{S'} = k(a_1x_1 - a_2y_2)^2$ *where* $k = b_1/a_1$;
- $r_{S'} = 16k(a_1x_1 + c_2y_2)^2$ *where* $k = b_1/a_1$.

Proof. Symbolic computation can verify the factorization once the coefficients are properly parametrized. For each of the systems, respectively set $\{b_1 = -a_1, c_2 = -a_2\}$, $\{b_1 = ka_1, c_2 = ka_2, f_1 = 0\}$, and $\{b_1 = ka_1, a_2 = kc_2, f_1 = 1\}$.

Now assume $r_{S'}$ is reducible. From Corollary 7.8 we have $\deg_{y_2}^{x_1}(r_{S'}) = (2,2)$. By regarding $y_2$ as an unknown, symbolic computation shows the discriminant of $r_{S'}$ is

$$\Delta = C(1 - 4a_1b_1x_1^2)[(4f_1^2(a_1a_2 - b_1c_2) + (f_1^2 - 1)^2(a_2b_1 - a_1c_2))x_1 - 2f_1(f_1^2 - 1)(a_2 + c_2)]^2$$

where $C$ is a positive constant. Since $S'$ has no constant branch, $r_{S'}$ is reducible if and only if $r_{S'}$ admits a factor of double degree $(1,1)$, i.e. $\Delta$ is a square in $\mathbb{C}[x_1]$. Notice that $(1 - 4a_1b_1x_1^2)$ always admits two distinct roots, so we must have

$$4f_1^2(a_1a_2 - b_1c_2) + (f_1^2 - 1)^2(a_2b_1 - a_1c_2) = 2f_1(f_1^2 - 1)(a_2 + c_2) = 0$$

to ensure $\Delta$ is a square, which leads to the three systems stated in the lemma. □

Lemma 7.11. *For a reduced coupling* $S' = (g_5^{(1)}, h^{(1)}, g_3^{(2)}, h^{(2)})$, $r_{S'}$ *is reducible if and only if one of the following systems holds.*

- $\{a_1c_1 = a_2b_2, f_1 = 0\}$;
- $\{a_1c_1a_2b_2 = \frac{1}{16}, f_1 = 1\}$.

*In particular, when* $\{a_1c_1 = a_2b_2, f_1 = 0\}$, *the factorization of* $r_{S'}$ *is as follow.*

- $r_{S'} = (a_2y_2 - a_1x_1)(x_1y_2 - k)$ *where* $k = b_2/a_1$.

Proof. Symbolic computation can verify the factorization once the coefficients are properly parametrized. For each of the systems, respectively set $\{c_1 = ka_2, b_2 = ka_1, f_1 = 0\}$ and $\{c_1 = \frac{-k^2}{4a_1}, b_2 = \frac{-1}{4k^2a_2}, f_1 = 1\}$.

Now assume $r_{S'}$ is reducible. Regarding $x_1, y_2$ as unknowns, we have respectively the discriminants of $G_5^{(1)}$ and $g_3^{(2)}$

$$\Delta_1 = 4(f_1^2 - a_1c_1(f_1^2-1)^2)x_2^2 + 4f_1(f_1^2-1)(1+4a_1c_1)x_2 + (f_1^2-1)^2 - 16a_1c_1f_1^2,$$
$$\Delta_2 = -4a_2b_2x_2^2 + 1.$$

On the other hand, Theorem 7.7 tells $\deg_{y_2}^{x_1}(r_{S'}) = (2,2)$. Since $S'$ has no constant branch, $r_{S'}$ has a factor $r$ which is in the form $r = (sx_1 + t)y_2 - (ux_1 + v)$. Notice that, as resultants, $\{G_5^{(1)}, r_{S'}\} \subset S'$, so $Z(S')$ has a component in the form

$$W = Z(g_5^{(1)}, h^{(1)}, g_3^{(2)}, h^{(2)}, G_5^{(1)}, r).$$

It is clear that $W$ is a (1-dimensional) curve. Choose $p \in W$ around which $W$ can be locally parametrized as $\gamma(x_2)$. In particular, there are rational functions $P_i, Q_i \in \mathbb{C}(x_2)$ for $i = 1,2$ such that $x_1 = P_1 + Q_1\sqrt{\Delta_1}$ and $y_2 = P_2 + Q_2\sqrt{\Delta_2}$. Moreover, given $r|_W \equiv 0$, $y_2 = \frac{ux_1+v}{sx_1+t}$ always holds. This implies $\sqrt{\Delta_1\Delta_2} \in \mathbb{C}[x_2]$ hence $\Delta_1\Delta_2$ must be a square in $\mathbb{C}[x_2]$. Consequently, $\Delta_1$ must be a multiple of $\Delta_2$ given the latter admits two distinct roots. This leads to

$$4f_1(f_1^2-1)(1+4a_1c_1) = 0, \quad \frac{4(f_1^2 - a_1c_1(1-f_1^2)^2)}{-4a_2b_2} = \frac{(1-f_1^2)^2 - 16a_1c_1f_1^2}{1}.$$

Finally, we have $1 + 4a_1c_1 \neq 0$ thanks to Corollary 2.5, so one of the systems stated in the lemma must hold. □

**Remark 2.** In Lemma 7.11, we did not display the factorization of $r_{S'}$ when $\{a_1c_1a_2b_2 = \frac{1}{16}, f_1 = 1\}$ because it is useless for mesh construction. To see this, we should first notice that $a_1c_1$ and $a_2b_2$ are both positive. The reason is $16a_1c_1a_2b_2 = 1$ infers $a_1c_1$ and $a_2b_2$ share the same sign and it is either $|4a_1c_1| \geq 1$ or $|4a_2b_2| \geq 1$. Given Corollary 2.5, we have $a_1c_1$ and $a_2b_2$ both positive. On top of this, $f_1 = 1$ implies $\Delta_1/\Delta_2 < 0$. So, in $Z(S')$, there are at most finitely many points at which $x_1$ and $y_2$ are both real. This means the corresponding mesh is not flexible in $\mathbb{R}^3$.

### 7.2 Singular flexible meshes with two isograms

This type comes from the first combination in Theorem 7.9, which naturally splits the corresponding matching into two cases, **adjacent-isogonal singular matching** ($j_2 = j_3 = 1$) and **opposite-isogonal singular matching** ($j_1 = j_3 = 1$). For the first case, it is more convenient to consider the relabeling in which $j_3 = j_4 = 1$ so that $\deg_{x_3}^{x_1}(R_{S_1'}) = (2,2)$ and $\deg_{x_3}^{x_1}(R_{S_2'}) = (1,1)$ by Corollary 7.8. So $R_{S_1'}$ must be reducible if $\gcd(R_{S_1'}, R_{S_2'}) \neq 1$.

**Example 7.1** (General case for adjacent-isogonal singular matching)**.** Choose one of the following four systems to set up $\{f_1, a_i, b_i, c_i\}$ for $i = 1,2$.

- $\{b_1 = -a_1, c_2 = -a_2, c_1 = b_2 = 0\}$;
- $\{a_1c_2 = a_2b_1, f_1 = c_1 = b_2 = 0\}$;
- $\{a_1a_2 = b_1c_2, f_1 = 1, c_1 = b_2 = 0\}$;
- $\{a_1c_1 = a_2b_2, f_1 = b_1 = c_2 = 0\}$.

According to Lemma 7.10 and 7.11, choose a factor of $r_{S'}$ and express it in the form $(px_1 + q)y_2 - (mx_1 + n)$ where $(m, n, p, q)$ are in terms of $(f_1, a_i, b_i, c_i)$. Set $e_i = 0$ for $i = 1,2$, $b_i = c_i = 0$ for $i = 3,4$. Set $f_2 \in (-1,1]$, and $(a_3, e_3, a_4)$ arbitrary but nonzero, and $k_3 = \frac{-1\pm\sqrt{1-4a_3e_3}}{2a_3}$. Let

$$N_1 = \begin{pmatrix} m & n \\ p & q \end{pmatrix},\ N_2 = \begin{pmatrix} 1-f_2^2 & -2f_2 \\ 2f_2 & 1-f_2^2 \end{pmatrix}.$$

Now, go through all the procedures in Example 5.1 starting from equation (13) till the end.

The idea of the above example is quite similar to Example 5.1: To construct Möbius transformations $N_i$ such that $N_4N_3N_2N_1(x_1) = x_1$, see proof of Theorem 5.1.

For opposite-isogonal singular matching, the general construction is derived from technical calculations. However, the idea is easy to follow. We strongly suggest the readers accept all the computational results till they capture the main idea, and then verify the calculation details efficiently through symbolic computation.

According to Theorem 7.9 we may assume $S_1' = (g_1^{(1)}, h^{(1)}, g_3^{(2)}, h^{(2)})$ and $S_2' = (g_1^{(3)}, h^{(3)}, g_5^{(4)}, h^{(4)})$. It is clear that

$$Z(S_1') \underset{\{x_1,x_3\}}{\times} Z(S_2') = Z(\tilde{S}_1') \underset{\{x_4,y_4\}}{\times} Z(g_5^{(4)}) \tag{18}$$

where $\tilde{S}_1' = (h^{(4)}, g_1^{(1)}, h^{(1)}, g_3^{(2)}, h^{(2)}, g_1^{(3)}, h^{(3)})$. Notice that, except for $g_3^{(2)}$, every generator of $\tilde{S}_1'$ induces a Möbius transformation between its arguments. Hence, in $Z(\tilde{S}_1')$, we have $x_2 = \frac{uy_4+v}{sy_4+t}$ and $y_2 = \frac{Ux_4+V}{Sx_4+T}$. Specifically,

$$\frac{uy_4+v}{sy_4+t} = \frac{(k_1F_4-F_1)y_4+(F_1F_4+k_1)}{(k_1F_1F_4+1)y_4+(k_1F_1-F_4)},\quad \frac{Ux_4+V}{Sx_4+T} = \frac{(F_2-k_3F_3)x_4+(F_2F_3+k_3)}{(k_3F_2F_3+1)x_4+(F_3-k_3F_2)}, \tag{19}$$

where $F_i = \frac{2f_i}{1-f_i^2}$ ($h^{(i)} = 0 \Leftrightarrow y_i = \frac{x_{i+1}+F_i}{1-F_ix_{i+1}}$). The corresponding Möbius map is

$$\mathcal{M}(g_3^{(2)}(x_2,y_2)) = (sy_4+t)(Sx_4+T)^2 g_3^{(2)}\left(\frac{uy_4+v}{sy_4+t}, \frac{Ux_4+V}{Sx_4+T}\right) \in \mathbb{C}[x_4,y_4].$$

We naturally have $Z(\mathcal{M}(g_3^{(2)})$ the $(x_4,y_4)-$projection of $Z(\tilde{S}_1')$. Given $g_5^{(4)}$ irreducible, the fiber product (18) being an infinite set requires $\mathcal{M}(g_3^{(2)})/g_5^{(4)}$ being a constant. Express $g_3^{(2)}, g_5^{(4)}$ in a matrix form:

$$g_3^{(2)} = \begin{pmatrix} y_2^2 \\ y_2 \\ 1 \end{pmatrix}^{\mathrm{T}} \begin{pmatrix} a_2 & 0 \\ 0 & 1 \\ b_2 & 0 \end{pmatrix} \begin{pmatrix} x_2 \\ 1 \end{pmatrix},\ g_5^{(4)} = \begin{pmatrix} x_4^2 \\ x_4 \\ 1 \end{pmatrix}^{\mathrm{T}} \begin{pmatrix} a_4 & 0 \\ 0 & 1 \\ c_4 & 0 \end{pmatrix} \begin{pmatrix} y_4 \\ 1 \end{pmatrix}.$$

Therefore, $\mathcal{M}(g_3^{(2)})$ can be rewritten as $\mathcal{M}(g_3^{(2)}) = (x_4^2, x_4, 1)M(y_4,1)^{\mathrm{T}}$ where

$$M = \begin{pmatrix} m_{11} & m_{12} \\ m_{21} & m_{22} \\ m_{31} & m_{32} \end{pmatrix} = \begin{pmatrix} U^2 & SU & S^2 \\ 2UV & SV+UT & 2ST \\ V^2 & VT & T^2 \end{pmatrix} \begin{pmatrix} a_2 & 0 \\ 0 & 1 \\ b_2 & 0 \end{pmatrix} \begin{pmatrix} u & v \\ s & t \end{pmatrix}. \tag{20}$$

Given $\mathcal{M}(g_3^{(2)})/g_5^{(4)} \in \mathbb{R}$, we need $m_{12} = m_{21} = m_{32} = 0$. In addition, since the matrices of Möbius transformations are considered in $\mathrm{PGL}(2,\mathbb{R})$, it is harmless to assume that

$$d := ut - sv,\ D := UT - SV,\ \{d, D\} \subset \{\pm 1\}.$$

Hence, by calculating $m_{ij}$ we have

$$\begin{cases} (S^2b_2 + U^2a_2)v + SUt = 0, \\ (T^2b_2 + V^2a_2)v + VTt = 0, \end{cases} \begin{cases} 2(STb_2 + UVa_2)u + (SV+TU)s = 0, \\ (tu - vs = d) \wedge (STb_2 = UVa_2), \end{cases} \tag{21}$$

a linear system of $(v, t)$ and a linear system of $(u, s)$. The last equation $STb_2 = UVa_2$ ensures nonzero solutions for $(v, t)$ so that $tu - vs = d$ is solvable for $(u, s)$. Generally, for preassigned $(U, V, S, T)$ such that $UT - SV = D$ and $STb_2 = UVa_2$, system (21) can provide feasible $(u, v, s, t)$ having $m_{12} = m_{21} = m_{32} = 0$. So $(a_4, c_4)$ can be obtained from

$$\frac{a_4}{m_{11}} = \frac{1}{m_{22}} = \frac{c_4}{m_{31}}. \tag{22}$$

**Example 7.2** (General case for opposite-isogonal singular matching)**.** Set $d, D \in \{\pm 1\}$, $(a_1, a_2, b_2, a_3)$ nonzero, and $(b_1, c_1, c_2, e_2, b_3, c_3, b_4, e_4)$ all 0. Pick $U, V, S, T$ such that $UT - SV = D$ and $STb_2 = UVa_2$. Therefore, system (21) is solvable for $(u, v, s, t)$: choose $n \neq 0$ and set

$$(v, t) = n(-SU, S^2b_2 + U^2a_2), \ (u, s) = \left(\frac{d(SV + TU)}{nD(U^2a_2 - S^2b_2)}, \frac{2d(STb_2 + UVa_2)}{nD(S^2b_2 - U^2a_2)}\right).$$

Let $(p, q)$ be such that

$$dp^2 + (u^2 + v^2 + s^2 + t^2)p + d = 0, \ Dq^2 + (U^2 + V^2 + S^2 + T^2)q + D = 0.$$

Notice that $p, q$ are real given $(u^2 + v^2 + s^2 + t^2)^2 \geq 4|ut - sv|^2 = 4$ etc.. To solve system (19), set

$$(k_1, k_3, F_1, F_2, F_3, F_4) = \left(p, q, -\frac{dp + u^2 + v^2}{su + tv}, \frac{Dq + U^2 + V^2}{SU + TV}, -\frac{Uq + T}{Vq - S}, \frac{up + t}{vp - s}\right).$$

Use $F_i = \frac{2f_i}{1-f_i^2}$ to recover $f_i \in (-1, 1]$, and for $i = 1, 3$, adjust the value of $a_i$ so that $e_i = -a_ik_i^2 - k_i \neq 0$. Finally, calculate matrix (20) and set $(a_4, c_4) = \left(\frac{m_{11}}{m_{22}}, \frac{m_{31}}{m_{22}}\right)$ from equation (22). Start over if inequalities (16) do not hold for all $i$.

### 7.3 Singular flexible meshes with four deltoids

Now we consider the second combination in Theorem 7.9. The matching corresponding to this type of mesh is said to be **deltoidal**. Specifically, the purpose of Lemma 7.10 and 7.11 is to divide this type into **reducible deltoidal matching** and **irreducible deltoidal matching**. Let $S_1', S_2'$ be reduced couplings of a deltoidal matching. According to Corollary 7.8, $\deg_{x_3}^{x_1}(R_{S_1'}) = \deg_{x_3}^{x_1}(R_{S_2'}) = (2, 2)$. Reducible deltoidal matching requires $R_{S_1'}, R_{S_2'}$ share a common factor of double degree $(1, 1)$, irreducible deltoidal matching requires $R_{S_1'}, R_{S_2'}$ both irreducible hence only differ by a constant.

**Example 7.3** (Typical case for reducible deltoidal matchings)**.** Proceed Example 7.1 and stop when $(m, n, p, q)$ are obtained! Set $e_i = 0$ for $1 \leq i \leq 4$. Now, the reader may pick any one of the following two options to construct the remaining coefficients.

Op. 1. Set $\{a_3, a_4\}$ nonzero and $\{b_3, c_4\}$ all 0. Let $\begin{pmatrix} u & v \\ s & t \end{pmatrix} = \begin{pmatrix} m & n \\ p & q \end{pmatrix}$, and go through all the procedures in Example 5.1 starting from equation (13) and stop when $(f_3, k_4, f_4)$ are obtained! Then update the values $f_2 \leftarrow f_3$ and $f_3 \leftarrow 0$. Finally, set $c_3 = k_4a_4$ and $b_4 = k_4a_3$.

Op. 2. Set $c_3 = b_4 = 0$, $f_2, f_4 \in (-1, 1]$. Let

$$\begin{pmatrix} u & v \\ s & t \end{pmatrix} = \begin{pmatrix} 1 - f_2^2 & -2f_2 \\ 2f_2 & 1 - f_2^2 \end{pmatrix} \begin{pmatrix} m & n \\ p & q \end{pmatrix} \begin{pmatrix} 1 - f_4^2 & -2f_4 \\ 2f_4 & 1 - f_4^2 \end{pmatrix}. \tag{23}$$

Adjust $f_2, f_4$, if necessary, such that there is a real number $z$ satisfying $sv^4z^4 + (4uvt - 2sv^2)z^2 + s = 0$. Set

$$f_3 = \frac{-1 + \sqrt{z^2v^2 + 1}}{zv}, \ a_3 = -b_3 = \frac{zt}{1 - z^2v^2}, \ a_4 = -c_4 = \frac{s(z^2v^2 - 1)}{4zvt}. \tag{24}$$

Start over if inequalities (6) do not hold for all $i$.

**Illustration:** A reducible deltoidal matching consists of reduced couplings $(S_1', S_2')$ such that $(r_{S_1'}, r_{S_2'})$ are reducible. Specifically, in $W := Z(S_1') \times_{\{x_1,x_3\}} Z(S_2')$ we have

$$h^{(4)}(y_4, x_1) = r_{S_2'}(x_3, y_4) = h^{(2)}(y_2, x_3) = r_{S_1'}(x_1, y_2) = 0. \tag{25}$$

Example 7.3 start with a reducible $r_{S_1'}$ having a factor $(px_1 + q)y_2 - (mx_1 + n)$. Observing the setups of the coefficients, Op. 1 considered $r_{S_2'}$ with a factor $(x_3y_4 - k_4)$ from Lemma 7.11, hence equation (25) provides consecutive Möbius transformations

$$x_1 = \begin{pmatrix} 1-f_4^2 & -2f_4 \\ 2f_4 & 1-f_4^2 \end{pmatrix} \begin{pmatrix} 0 & k_4 \\ 1 & 0 \end{pmatrix} \begin{pmatrix} 1-f_2^2 & -2f_2 \\ 2f_2 & 1-f_2^2 \end{pmatrix} \begin{pmatrix} m & n \\ p & q \end{pmatrix} (x_1).$$

In order to make $W$ an infinite set, the above matrix production should result in a scalar matrix. The construction for $(f_2, f_4, k_4)$ goes the same way as in Example 5.1 for $(f_3, f_4, k_4)$. On the other hand, Op. 2 considered $r_{S_2'}$ with a factor from the first case of Lemma 7.10. Notice that $h^{(4)} = h^{(2)} = 0$ simplifies equation (25) as

$$r_{S_2'}(x_3, y_4) = r_{S_1'}\left(\frac{(1-f_4^2)y_4 - 2f_4}{2f_4y_4 + 1 - f_4^2}, \frac{(1-f_2^2)x_3 + 2f_2}{1 - f_2^2 - 2f_2x_3}\right) = \frac{f(x_3, y_4)}{g(x_3, y_4)} = 0,$$

where $f, g \in \mathbb{C}[x_3, y_4]$. A series of Möbius transformations in (23) convert the factor $(px_1 + q)y_2 - (mx_1 + n)$ of $r_{S_1'}$ to the factor $(sy_4 + t)x_3 - (uy_4 + v)$ of $f(x_3, y_4)$. Therefore, $r_{S_2'}$ must also take $(sy_4 + t)x_3 - (uy_4 + v)$ as a factor so that $W$ can be an infinite set. Symbolic computation can demonstrate the correctness of Op. 2 by factorizing $r_{S_2'}$ with the coefficients settled by equation (24).

According to Lemma 7.10 and 7.11, the only case of the factor of $r_{S_2'}$ that is not covered is $(x_3 + ky_4)$. However, by setting $\bar{y}_4 = -\frac{1}{y_4}$, it is converted to Op. 1 since $(x_3 + ky_4) = 0 \Leftrightarrow (x_3\bar{y}_4 - k) = 0$. And this is why we call Example 7.3 'typical case'.

Proposition 7.12. *Let $S_1' = (g_3^{(1)}, h^{(1)}, g_5^{(2)}, h^{(2)})$ and $S_2' = (g_3^{(3)}, h^{(3)}, g_5^{(4)}, h^{(4)})$ be two reduced couplings correspond to a matching $M$ such that $Z(S_1') \times_{\{x_1,x_3\}} Z(S_2')$ is an infinite set. Then $M$ must be reducible.*

Proof. Take $x_3$ as an unknown, symbolic computation shows the ratio of the discriminants of $R_{S_1'}$ and $R_{S_2'}$ is

$$\frac{\Delta_1}{\Delta_2} = \frac{(1 - 4a_1b_1x_1^2)p^2}{4(f_4^2 - a_4c_4(1-f_4^2)^2)x_1^2 + 4f_4(f_4^2 - 1)(1 + 4a_4c_4)x_2 + (1-f_4^2)^2 - 16a_4c_4f_4^2}$$

where $p \in \mathbb{C}(x_1)$ is a rational function of $x_1$. Given $Z(S_1') \times_{\{x_1,x_3\}} Z(S_2')$ is an infinite set, we have $\gcd(R_{S_1'}, R_{S_2'}) \neq 1$. Suppose $R_{S_1'}$ is irreducible, $R_{S_2'}$ must be a multiple of $R_{S_1'}$, which implies $\frac{\Delta_1}{\Delta_2} \in \mathbb{C}$. Considering $1 - 4a_1b_1x_1^2$ has two distinct roots, it must be $p \in \mathbb{C}$ and

$$4f_4(f_4^2 - 1)(1 + 4a_4c_4) = 0, \quad \frac{4(f_4^2 - a_4c_4(1-f_4^2)^2)}{-4a_1b_1} = \frac{(1-f_4^2)^2 - 16a_4c_4f_4^2}{1}. \tag{26}$$

Set $\tilde{S}_1' = (g_5^{(4)}, h^{(4)}, g_3^{(1)}, h^{(1)})$ and $\tilde{S}_2' = (g_5^{(2)}, h^{(2)}, g_3^{(3)}, h^{(3)})$. Thus

$$Z(S_1') \underset{\{x_1,x_3\}}{\times} Z(S_2') = Z(\tilde{S}_1') \underset{\{x_2,x_4\}}{\times} Z(\tilde{S}_2').$$

Review the end of the proof of Lemma 7.11, equation (26) says $r_{\tilde{S}_1'}$ (hence $R_{\tilde{S}_1'}$) is reducible. This means $M$ is always reducible from the perspective of relabeling. □

Let $S_1', S_2'$ be two reduced couplings corresponding to an irreducible deltoidal matching. The above proposition guides us to the only possible situation:

$$S_1' = (g_5^{(1)}, h^{(1)}, g_3^{(2)}, h^{(2)}),\ \ S_2' = (g_3^{(3)}, h^{(3)}, g_5^{(4)}, h^{(4)}).$$

The construction idea is simple. Set up $S_1'$ such that $R_{S_1'}$ is irreducible (by avoiding the conditions in Lemma 7.11), construct $S_2'$ such that $R_{S_2'}$ and $R_{S_1'}$ only differ by multiplying a constant. Once the coefficients of $S_1'$ are evaluated in $\mathbb{R}$, using numerical methods to obtain $S_2'$ is effortless since $R_{S_2'}$ can be given explicitly by symbolic computation right away. However, to obtain a general solution, like what we did to all previous cases, is much more difficult. In fact, the coefficients must satisfy

$$\begin{cases} \frac{a_1}{c_1} = \frac{4(a_4c_4(f_4^2-1)^2-f_4^2)}{16a_4c_4f_4^2-(f_4^2-1)^2}, & \frac{a_2}{b_2} = \frac{4(a_3b_3(f_2^2-1)^2-f_2^2)}{16a_3b_3f_2^2-(f_2^2-1)^2}, \\ \frac{f_1}{1-f_1^2} \in \left\{\pm\frac{\sqrt{a_1c_1}-\sqrt{a_2b_2}}{4\sqrt{a_1c_1a_2b_2}+1}, \pm\frac{\sqrt{a_1c_1}+\sqrt{a_2b_2}}{4\sqrt{a_1c_1a_2b_2}-1}\right\}, \end{cases} \tag{27}$$

which requires $a_1c_1 > 0$ and $a_2b_2 > 0$ to obtain an $f_1 \in (-1, 1]$. The above result comes from an extremely complicated calculation which will not be presented here. However, there is no harm to use it as a hint for constructions as long as the validity of the output can be easily tested. Hence, in this article, we only provide a special case of irreducible deltoidal matching to demonstrate its existence.

**Example 7.4** (Special case for irreducible deltoidal matchings)**.** Set

$$\{b_1, e_1, c_2, e_2, f_3, c_3, e_3, b_4, e_4\} = \{0\},\ \ f_2 = f_4 = \sqrt{2} - 1,$$

and $\{a_1, a_2, a_3\}$ arbitrary but $a_1a_2a_3^2 < 0$. According to equation (27), set $c_1 = a_1, b_2 = a_2$, and $f_1 \in (-1, 1]$ such that $\frac{f_1}{1-f_1^2} = \frac{a_1+a_2}{4a_1a_2-1}$. Adjust the values of $(a_1, a_2)$ to avoid the conditions in Lemma 7.11 so that $r_{S_1'}$ is indeed irreducible. Finally, set

$$b_3 = \frac{(a_1 + a_2)(1 - 4a_1a_2)}{4a_1a_3(4a_2^2 + 1)},\ \ a_4 = \pm a_3\sqrt{-\frac{a_1(4a_2^2 + 1)}{a_2(4a_1^2 + 1)}},\ \ c_4 = -\frac{b_3a_4}{a_3}.$$

Start over if inequalities (6) do not hold for all $i$.

**Remark 3.** Let $S_1' = (g_1^{(1)}, h^{(1)}, g_3^{(2)}, h^{(2)})$ and $S_2' = (g_1^{(3)}, h^{(3)}, g_5^{(4)}, h^{(4)})$ be reduced couplings constructed in Example 7.4. Symbolic computation can efficiently test its correctness by showing $R_{S_1'}/R_{S_2'}$ is a constant. In addition, let $\tilde{S}_1' = (g_5^{(4)}, h^{(4)}, g_3^{(1)}, h^{(1)})$, Corollary 7.8 says $\deg_{x_2}^{x_4}(R_{\tilde{S}_1'}) = (4, 1)$. Since $R_{\tilde{S}_1'}$ has no factor in $\mathbb{C}[x_4]$, it is naturally irreducible. This means the matching is indeed irreducible no matter how we label the polynomials.

## 8 SUMMARY AND FUTURE WORK

This article thoroughly studied and classified flexible meshes with reducible quadrilaterals. The classification contains following typical classes.

$$\begin{cases} \text{(Non-singular) isogonal,} \\ \text{(Singular) constant,} \\ \text{(Singular) non-constant} \begin{cases} \text{isogonal} \begin{cases} \text{adjacent,} \\ \text{opposite,} \end{cases} \\ \text{detoidal} \begin{cases} \text{reducible,} \\ \text{irreducible.} \end{cases} \end{cases} \end{cases} \tag{28}$$

Any other case can be converted to one of the above by the tricks in Section 4. Moreover, systematic constructions are given to recover all possible meshes of every class in diagram (28) except for irreducible detoidal meshes due to high computational complexity. However, we still discovered a special case for such a class that can provide plentiful flexible samples for readers to explore.

This work is the first piece of the blueprint of the classification problem on 3-by-3 flexible quadrilateral meshes with skew faces. It is also a milestone considering its history and difficulty: The first mesh of such type was discovered in [6], decades after the problem had been posed. The future work is split into two parts.

1. Flexible meshes with irreducible quadrilaterals;
2. Flexible meshes with hybrid (reducible and irreducible) quadrilaterals.

## A DEFINITION OF THE ANGLES

Let us consider a planar mesh that is partially depicted in Fig.5 left. We relocate the starting points of all vectors to the center of a unit sphere to obtain a spherical linkage in Fig.5 right. We take the common orientation of the sphere with surface normals pointing outwards. In quad $Q_1 : (\lambda_1, \delta_1, \mu_1, \gamma_1)$, the dihedral angles $\alpha_1, \beta_1 \in [0, 2\pi)$ are well defined by right-hand rule: the anticlockwise rotating angle from $\lambda_1$ (or $\delta_1$) to $\gamma_1$ (or $\lambda_1$). $\alpha_2, \beta_2 \in [0, 2\pi)$ are similarly defined.

For skew meshes (Fig.2), we formalize the angle determination as follows.

Let $\boldsymbol{v_i}$ be the surface normal at the common vertex of $Q_i$ and $Q_{i+1}$. Take Fig.3 as an example, a neighborhood of $Q_1 \cap Q_2$ can be stereographically projected (angle preserving) to its tangent plane (with the same orientation, see Fig.6). Let $T(\boldsymbol{v_i})$ be the oriented tangent plane with surface normal $\boldsymbol{v_i}$. Set $-\lambda_i$ and $-\delta_i$, respectively, are the reflections of $\lambda_i$ and $\delta_i$ about $\boldsymbol{v_i}$. Introduce the operator

$$< \cdot, \cdot >_{\boldsymbol{v_i}} : T(\boldsymbol{v_i})^2 - \{\boldsymbol{0}\} \to [0, 2\pi);$$

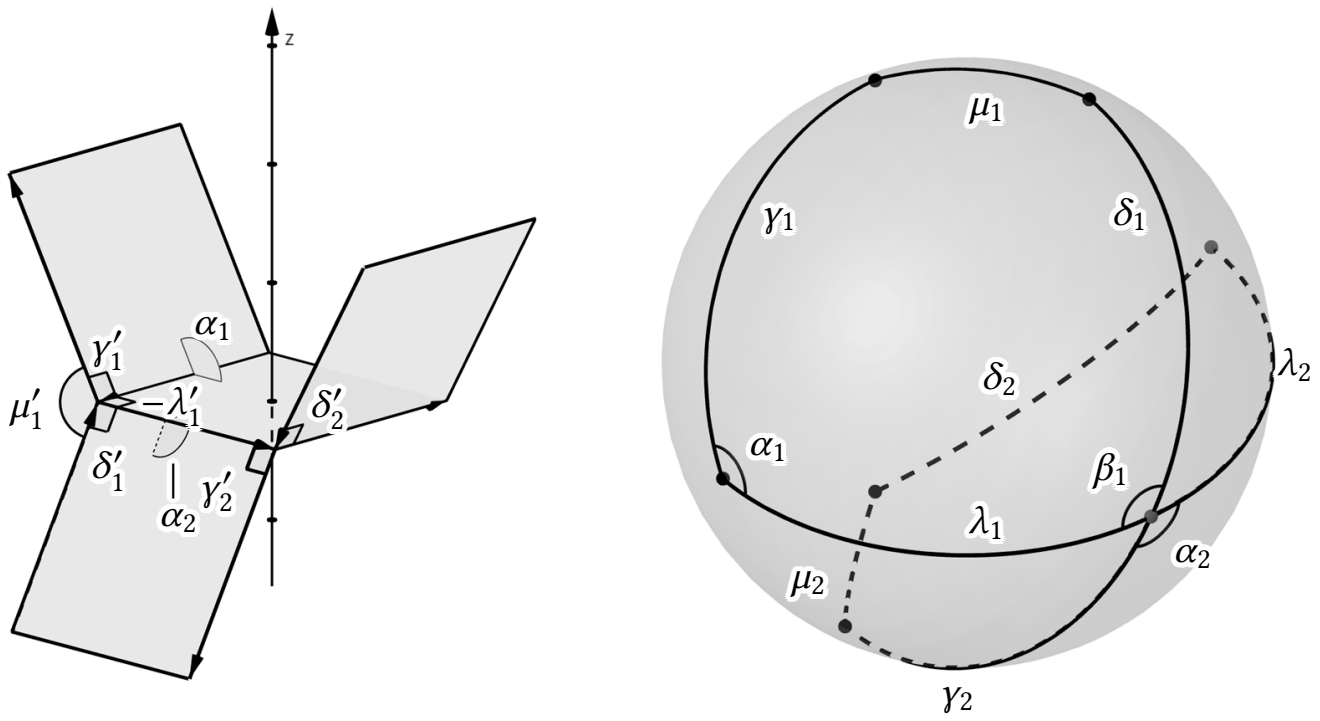

Fig. 5. Notice that $(\lambda_i, \gamma_i, \mu_i, \delta_i)$ and $(\lambda_i', \gamma_i', \mu_i', \delta_i')$ are complementary to $\pi$ respectively.

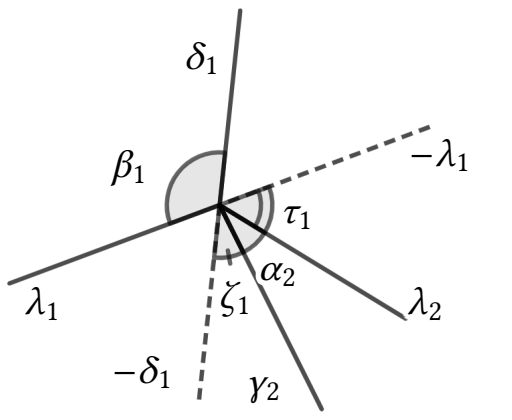

Fig. 6. Stereographic projection of the neighborhood near the '$*$' in Fig.3 left. $\boldsymbol{v_1}$ at the intersecting point is a vector perpendicular to the paper and pointing towards the reader. $-\lambda_1$ and $-\delta_1$, respectively, are the reflections of $\lambda_1$ and $\delta_1$ about $\boldsymbol{v_1}$.

$$< \boldsymbol{x}, \boldsymbol{y} >_{\boldsymbol{v_i}} := \text{The counterclockwise rotating angle from } \boldsymbol{x} \text{ to } \boldsymbol{y}.$$

$$\alpha_i := \begin{cases} < \lambda_i, \gamma_i >_{\boldsymbol{v_{i-1}}}, & i \in \{1,3\}, \\ < \gamma_i, \lambda_i >_{\boldsymbol{v_{i-1}}}, & i \in \{2,4\}, \end{cases} \qquad \beta_i := \begin{cases} < \delta_i, \lambda_i >_{\boldsymbol{v_i}}, & i \in \{1,3\}, \\ < \lambda_i, \delta_i >_{\boldsymbol{v_i}}, & i \in \{2,4\}, \end{cases}$$

where $\boldsymbol{v_0} = \boldsymbol{v_4}$. Thereafter, by setting

$$\tau_i := \begin{cases} < \lambda_{i+1}, -\lambda_i >_{\boldsymbol{v_i}}, & i \in \{1,3\}, \\ < -\lambda_i, \lambda_{i+1} >_{\boldsymbol{v_i}}, & i \in \{2,4\}, \end{cases} \qquad \zeta_i := \begin{cases} < -\delta_i, \gamma_{i+1} >_{\boldsymbol{v_i}}, & i \in \{1,3\}, \\ < \gamma_{i+1}, -\delta_i >_{\boldsymbol{v_i}}, & i \in \{2,4\}, \end{cases}$$

we have $\beta_i \equiv \alpha_{i+1} + \tau_i + \zeta_i \mod 2\pi$ for all $i \in \{1,2,3,4\}$. For instance,

$$\begin{aligned} \beta_2 &= < \lambda_2, \delta_2 >_{\boldsymbol{v_2}} = < -\lambda_2, -\delta_2 >_{\boldsymbol{v_2}} \\ &\equiv < -\lambda_2, \lambda_3 >_{\boldsymbol{v_2}} + < \lambda_3, \gamma_3 >_{\boldsymbol{v_2}} + < \gamma_3, -\delta_2 >_{\boldsymbol{v_2}} = \tau_2 + \alpha_3 + \zeta_2 \mod 2\pi. \end{aligned}$$